\documentclass[12pt,a4paper,twoside]{article}

\usepackage{graphicx}
\usepackage{amsfonts,amsmath,amssymb,amsthm,latexsym}
\usepackage{times,mathptm}
\usepackage{handelingen}

\newtheorem{proposition}{Proposition}[section]
\newtheorem{lemma}[proposition]{Lemma}

\newtheorem{theorem}[proposition]{Theorem}

\theoremstyle{definition}
\newtheorem{definition}[proposition]{Definition}

\theoremstyle{remark}
\newtheorem{remark}[proposition]{Remark}

\newcommand{\thlabel}[1]{\label{th:#1}}
\newcommand{\thref}[1]{Theorem~\ref{th:#1}}
\newcommand{\selabel}[1]{\label{se:#1}}
\newcommand{\seref}[1]{Section~\ref{se:#1}}
\newcommand{\lelabel}[1]{\label{le:#1}}

\newcommand{\prlabel}[1]{\label{pr:#1}}
\newcommand{\prref}[1]{Proposition~\ref{pr:#1}}

\newcommand{\delabel}[1]{\label{de:#1}}

\newcommand{\eqlabel}[1]{\label{eq:#1}}
\newcommand{\equref}[1]{(\ref{eq:#1})}

\renewcommand{\thefootnote}{\fnsymbol{footnote}}
\title{\MakeUppercase{On Reducibility of Mapping Class Group Representations: The $SU(N)$ Case}}
\author{J\o{}rgen~Ellegaard~Andersen$^1$ and Jens~Fjelstad$^2$}
\date{}
\renewcommand{\date}{\vspace{-5mm}}
\begin{document}
\maketitle \vspace*{-3mm}\relax
\renewcommand{\thefootnote}{\arabic{footnote}}

\noindent \textit{\small CTQM, Department of Mathematical Sciences, Aarhus University, Ny Munkegade, bldg. 1.530,
8000 Aarhus C, Denmark}\\
\noindent \textit{\small $^1$ e-mail: andersen@imf.au.dk}\\
\noindent \textit{\small $^2$ e-mail: fjelstad@imf.au.dk}

\begin{abstract}
\noindent We review and extend the results of~\cite{AF} that give a condition for reducibility of quantum representations of mapping class groups constructed from Reshetikhin-Turaev type topological quantum field theories based on modular categories. This criterion is derived using methods developed to describe rational conformal field theories, making use of Frobenius algebras and their representations in modular categories. Given a modular category $\mathcal{C}$, a rational conformal field theory can be constructed from a Frobenius algebra $A$ in $\mathcal{C}$. We show that if $\mathcal{C}$ contains a symmetric special Frobenius algebra $A$ such that the torus partition function $Z(A)$ of the corresponding conformal field theory is non-trivial, implying reducibility of the genus $1$ representation of the modular group, then the representation of the genus $g$ mapping class group constructed from $\mathcal{C}$ is reducible for every $g\geq 1$. We also extend the number of examples where we can show reducibility significantly by establishing the existence of algebras with the required properties using methods developed by Fuchs, Runkel and Schweigert. As a result we show that the quantum representations are reducible in the $SU(N)$ case, $N>2$, for all levels $k\in\mathbb{N}$. The $SU(2)$ case was treated explicitly in~\cite{AF}, showing reducibility for even levels $k\geq 4$.
\end{abstract}

\section*{Introduction}
An important feature of three-dimensional topological quantum field theory (TQFT) is that it produces projective representations of mapping class groups (MCG's) of closed surfaces, possibly with a finite number of marked points. This paper is concerned with closed surfaces with no marked points, and representations of the corresponding MCG's given by Reshetikhin-Turaev type TQFT's~\cite{RT1,RT2,T}, so-called ``quantum representations''. A TQFT of this type is defined by a modular category. Denote the TQFT constructed from the category $\mathcal{C}$ by $\mathbf{tqft}_\mathcal{C}$, and the corresponding genus $g$ quantum representation by $V_g^\mathcal{C}$. Although there is an extensive body of literature on TQFT, a comparatively small amount is devoted to the properties of the associated quantum representations. Consequently there is little known in general, and what is known is restricted to a few special cases.
For the modular category constructed from representations of $U_q(\mathfrak{su}(2))$ when $q$ is a primitive $(k+2)$'nd root of unity, or alternatively the category of integrable highest weight representations of the untwisted affine Lie algebra $\widehat{\mathfrak{su}}(2)$ at level $k$, Roberts have shown~\cite{R} that the representations $V_g^\mathcal{C}$ are irreducible for all $g\geq 1$ when $k+2$ is prime. This result played a key role in proving that mapping class groups do not have Kazhadan's property $T$~\cite{A1}. More recently, Chen and Kerler have studied the case when $\mathcal{C}$ is the category of integrable highest weight representations of the affine Lie algebra $\widehat{\mathfrak{su}}(N)$ at level $k$ (or alternatively modular categories constructed from representations of $U_q(\mathfrak{su}(N))$ at a $(k+N)$'th root of unity), in particular for $g=1$~\cite{CK}. Apart from some explicit decompositions for $N=3$ and $g=1$, their results imply reducibility for any $g\geq 1$ when $N>2$, $k>2N$, and $\gcd(k,N)=1$.
For $g\geq 2$ there are remarkably few results. It has been shown that the $\mathfrak{su}(N)$ quantum representations are asymptotically faitful~\cite{A2}, i.e. the intersection of the kernels of the family of representations indexed by the level $k\in \mathbb{N}$ is trivial. In the geometric picture of~\cite{ADW, H} a decomposition of $\mathfrak{su}(2)$ quantum representations in subrepresentations for even level $k>2$ is known~\cite{AM}. For genus $g\geq 2$ it has been shown in the $\mathfrak{su}(2)$ case with the exception of levels $k=1, 2, 4, 8$ that the images of quantum representations are always infinite~\cite{M}.
Much more can be said about the genus $1$ representations, however. Partly because they are exceptional in more than one way, and partly due to having been under close scrutiny in connection with rational conformal field theory (RCFT).
Note that the genus $1$ mapping class group is $SL(2,\mathbb{Z})$. One of the exceptional features at genus $1$ is that any quantum representation factors through a finite group, $SL(2,\mathbb{Z}/N\mathbb{Z})$ for some integer $N$~\cite{CG,B}. As we have already seen, this does not happen generically for $g\geq 2$. Furthermore, genus $1$ representations, although projective, can always be lifted to linear representations as 
can be shown explicitly for any $\mathbb{C}$-linear modular category~\cite{BK} (see also \seref{tqft} below).
This particular property actually generalises to genus $2$, but for $g\geq 3$ it is known that $H^2(MCG_g,\mathbb{Q})\cong\mathbb{Q}$~\cite{KS}, so the extention class could
in general be non-zero, as it indeed is for the $SU(N)$ theories discussed in this paper.

There is a sizeable amount of literature on rational conformal field theory (RCFT) devoted to the commutant of $V_1^\mathcal{C}$ when $\mathcal{C}\cong\mathrm{Rep}(\mathcal{V})$ for a rational conformal vertex algebra $\mathcal{V}$ satisfying certain finiteness conditions. One problem relevant to RCFT is to find all elements of the commutant satisfying certain positivity and integrality constraints. More precisely, in a given basis $\{v_i\}_{i\in I}$ of the underlying vector space, the problem is to find all $|I|\times|I|$ matrices $Z$ commuting with matrices $S$ and $T$ generating the $SL(2,\mathbb{Z})$ representation such that all elements $Z_{i,j}$ are semi positive integers and with $Z_{0,0}=1$ where $0$ is a distinguished element of $I$. A solution to this problem is called a modular invariant, and finding all modular invariants for a given category $\mathcal{C}$ can be seen as an auxiliary task in finding all possible RCFT's associated to a given modular category (in particular, to a given vertex algebra). One immediate result is that if $\mathcal{C}$ contains simple objects that are not self-dual, then the genus $1$ representation is automatically reducible since one can then construct a non-trivial modular invariant. The commutant has been explicitly determined for $\mathcal{C}$ the category of integrable highest weight representations of $\widehat{\mathfrak{su}}(2)$ at level $k\in\mathbb{N}$, and for the category of highest weight representations of the Virasoro algebra with central charge $c\in(0,1)$~\cite{CIZ, GQ}, in both cases amounting to ADE classifications of modular invariants. In the first case, the results imply reducibility of the (genus $1$) quantum representations for every non-prime level $k$ (this follows from the determination of the whole commutant, non-trivial modular invariants only appear for even levels larger than $2$). Modular invariants are also completely classified for the case of $\widehat{\mathfrak{su}}(3)$ at level $k\in\mathbb{N}$~\cite{G}, implying reducibility of the modular group representations for all levels. Various other results are known, mainly related to representations of untwisted affine Lie algebras, but without complete classifications.

We begin in \seref{tqft} by establishing conventions and notation, and briefly reviewing some aspects of TQFT and quantum representations. In \seref{redcrit} we review relevant concepts and results from rational conformal field theory in the TQFT formalism and state the main result of~\cite{AF}, \thref{red}, a condition for reducibility of quantum representations. We also review the idea behind the proof, referring to~\cite{AF} for the full details. In \seref{sun} we first review the consequences of our reducibility result for the $\mathfrak{su}(2)$ case, and then continue by proving that the conditions of \thref{red} hold in the $\mathfrak{su}(N)$ case for $N\geq 3$ and for any level $k\in \mathbb{N}$, thus showing reducibility of the corresponding quantum representations (stated in \thref{redSUN}). The techniques we employ were developed in the series of papers~\cite{FRS1, FRS2, FRS3, FRS4, FjFRS1, FjFRS2, FrFRS1, FrFRS2}, and the proof presented here for the $\mathfrak{su}(N)$ case relies in particular on results from~\cite{FRS3} concerning Frobenius algebras constructed from direct sums of invertible objects in modular categories.

\section{TQFT and Quantum Representations}
\begin{figure}
	\begin{picture}(280,130)
		\put(0,95)				{$b_U=$}
		\put(30,80)			{\includegraphics[scale=0.3]{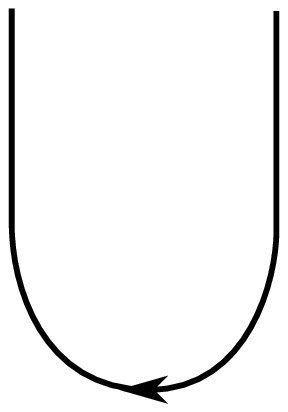}}
		\put(28,118)			{$U$}
		\put(50,118)			{$U^\vee$}
		\put(80,95)			{$d_U=$}
		\put(110,80)			{\includegraphics[scale=0.3]{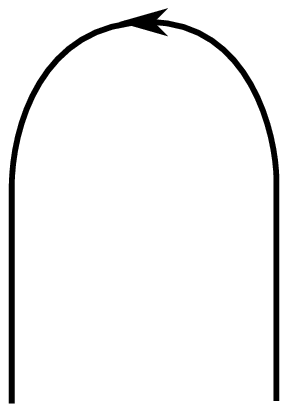}}
		\put(108,70)			{$U^\vee$}
		\put(130,70)			{$U$}
		\put(0,25)				{$\tilde b_U=$}
		\put(30,10)			{\includegraphics[scale=0.3]{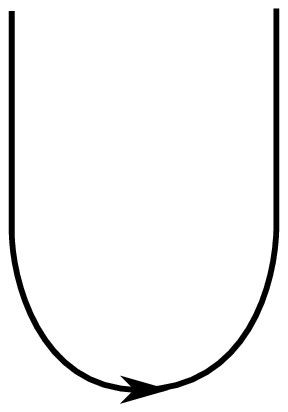}}
		\put(28,48)			{$U^\vee$}
		\put(50,48)			{$U$}
		\put(80,25)			{$\tilde d_U=$}
		\put(110,10)			{\includegraphics[scale=0.3]{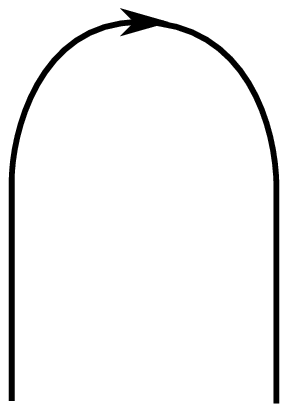}}
		\put(108,0)				{$U$}
		\put(130,0)			{$U^\vee$}
		\put(160,95)			{$c_{U,V}=$}
		\put(190,80)			{\includegraphics[scale=0.3]{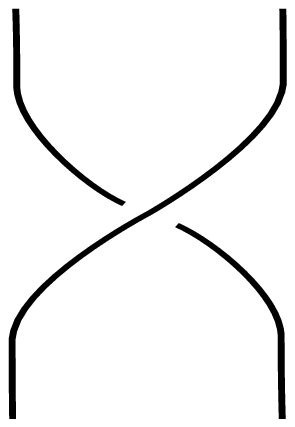}}
		\put(188,70)			{$U$}
		\put(210,70)			{$V$}
		\put(188,118)			{$V$}
		\put(210,118)			{$U$}
		\put(160,25)			{$c^{-1}_{U,V}=$}
		\put(190,10)			{\includegraphics[scale=0.3]{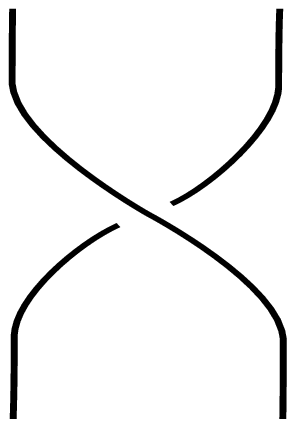}}
		\put(188,0)			{$V$}
		\put(210,0)			{$U$}
		\put(188,48)			{$U$}
		\put(210,48)			{$V$}
		\put(240,95)			{$\theta_U=$}
		\put(270,80)			{\includegraphics[scale=0.3]{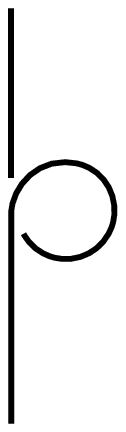}}
		\put(267,70)			{$U$}
		\put(267,118)			{$U$}
		\put(240,25)			{$\theta^{-1}_U=$}
		\put(270,10)			{\includegraphics[scale=0.3]{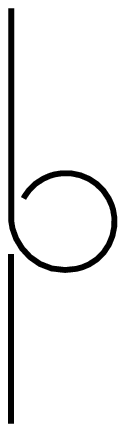}}
		\put(267,0)			{$U$}
		\put(267,48)			{$U$}
		\put(0,65)				{\line(1,0){280}}
		\put(70,0)				{\line(0,1){130}}
		\put(150,0)			{\line(0,1){130}}
		\put(230,0)			{\line(0,1){130}}
	\end{picture}
\caption{Structural morphisms for ribbon categories. $U,V\in\mathrm{Obj}(\mathcal{C})$}
\label{fig:strmorph}
\end{figure}

\selabel{tqft}
 We are interested in Reshetikhin-Turaev type topological quantum field theories (TQFT's) which are constructed from modular categories. Let us begin with conventions and notation.
 \paragraph{Modular Categories}
We will only deal with categories over $\mathbb{C}$, so by a modular category we mean an abelian, $\mathbb{C}$-linear, semisimple, ribbon category with simple tensor unit $\mathbf{1}$ and a finite number of isomorphism classes of simple objects, such that the braiding is maximally non-degenerate (to be explained below). We use graphical calculus for ribbon categories (see for instance~\cite{FRS1} for more elaborate discussion), where a string diagram is to be read from bottom to top. The structural morphisms (left \& right duality, braiding, twist) are displayed in figure~\ref{fig:strmorph}.
 
 We use the notation $\{U_i\}_{i\in I}$ for a chosen set of representatives of the isomorphism classes of simple objects in $\mathcal{C}$, where $U_0\cong\mathbf{1}$. The duality defines an involution $i\mapsto \bar\imath$ of the index set $I$ through $U_{\bar\imath}\cong U_i^\vee$.
The trace of an endomorphism $f\in\mathrm{Hom}(U,U)$ is defined as $\mathrm{tr}(f)=\tilde d_U\circ(f\otimes\mathrm{id}_{U^\vee})\circ b_U$. The modularity condition, i.e. maximal non-degeneracy of the braiding, is the requirement that the $|I|\times|I|$-matrix $s$ with elements 
 
\begin{center}
 \begin{picture}(220,80)
 	\put(0,40)				{ $s_{i,j}=\mathrm{tr}(c_{U_j,U_i}\circ c_{U_i,U_j})=$}
	\put(125,0)			{\includegraphics[scale=0.4]{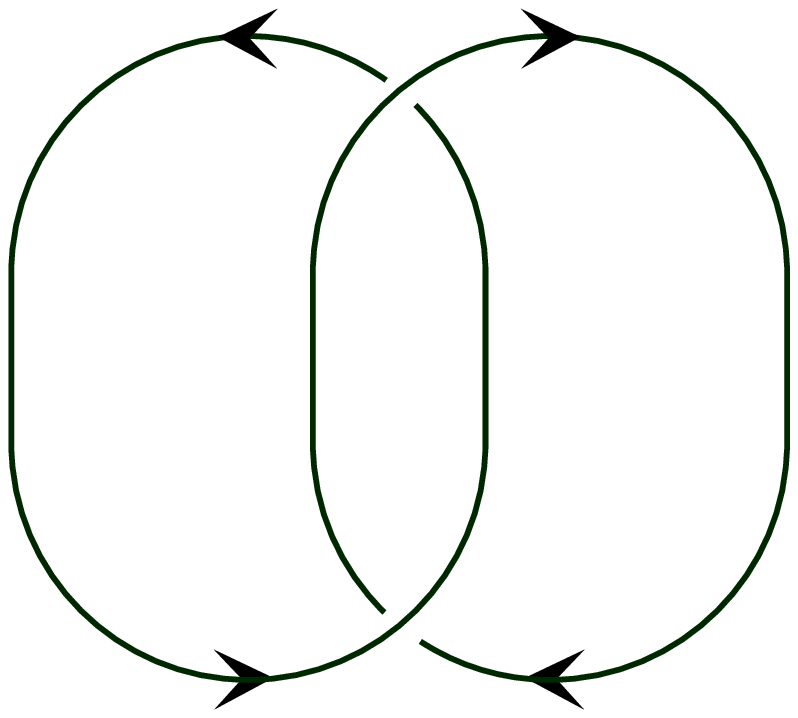}}
	\put(182,30)			{$U_i$}
	\put(148,30)			{$U_j$}
\end{picture}
\end{center}
 is invertible. 

By quantum dimension (or just dimension) of an object $U$ we mean the trace of the identity morphism of $U$, $\mathrm{dim}(U)=\mathrm{tr}(\mathrm{id}_U)$. We make the additional assumption that the dimension of any object is real and positive. The twist for a simple object satisfies $\theta_{U_i}=\theta_i\mathrm{id}_{U_i}$ for some $\theta_i\in\mathbb{C}^\times$. 

One shows~\cite{BK} that the $|I|\times|I|$ matrices $S$ and $T$ with elements $S_{i,j}=s_{i,j}/s_{0,0}$, $T_{i,j}=\zeta_\mathcal{C}\theta_i\delta_{i,j}$ with
$$\zeta_\mathcal{C}=\left(\frac{\sum_{i\in I}\theta_i(\mathrm{dim}(U_i))^2}{\sum_{i\in I}\theta^{-1}_i(\mathrm{dim}(U_i))^2}\right)^{1/6}$$
satisfy $(ST)^3=S^2$ and $S^4=1$. These are defining relations for $SL(2,\mathbb{Z})$ and thus define a matrix representation of the modular group, hence the name modular category.
It can also be shown that $S^2=C$, the charge conjugation matrix with elements $C_{i,j}=\delta_{\bar\imath,j}$.

\begin{figure}
	\begin{picture}(300,80)
		\put(0,10)				{\includegraphics[scale=0.5]{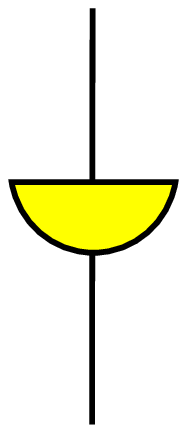}}
		\put(9,0)				{$V$}
		\put(9,73)				{$U$}
		\put(10,38)			{$e$}
		\put(50,10)			{\includegraphics[scale=0.5]{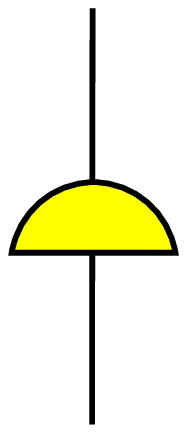}}
		\put(59,0)				{$U$}
		\put(59,73)			{$V$}
		\put(60,38)			{$r$}
		\put(120,10)			{\includegraphics[scale=0.5]{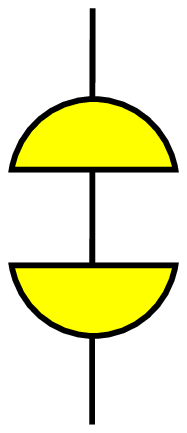}}
		\put(129,0)			{$V$}
		\put(135,38)			{$U$}
		\put(129,73)			{$V$}
		\put(130,26)			{$e$}
		\put(130,50)			{$r$}
		\put(155,40)			{$=$}
		\put(180,10)			{\includegraphics[scale=0.5]{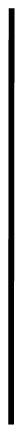}}
		\put(178,0)			{$V$}
		\put(178,73)			{$V$}
		\put(220,10)			{\includegraphics[scale=0.5]{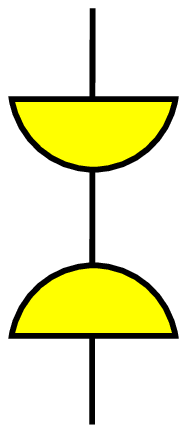}}
		\put(229,0)			{$U$}
		\put(235,38)			{$V$}
		\put(229,73)			{$U$}
		\put(230,26)			{$r$}
		\put(230,50)			{$e$}
		\put(255,40)			{$=$}
		\put(275,10)			{\includegraphics[scale=0.5]{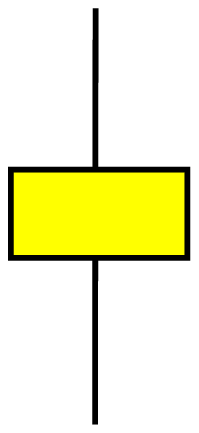}}
		\put(285,40)			{$p$}
		\put(284,0)			{$U$}
		\put(284,73)			{$U$}
	\end{picture}
\caption{Retract morphisms $e$ and $r$ for a retract $V\prec U$ corresponding to the split idempotent $p\in\mathrm{End}(U)$, and their relations.}
\label{fig:retract}
\end{figure}

The statement ``$V$ is a subobject of $U$'' is abbreviated $V\prec U$. Since $\mathcal{C}$ is abelian, every idempotent $p\in\mathrm{End}(U)$ (i.e. $p\circ p=p$) is split. We denote a corresponding retract by a triple $(V,e,r)$ where $V\in\mathrm{Obj}(\mathcal{C})$, $e\in\mathrm{Hom}(V,U)$ is a monic, $r\in\mathrm{Hom}(U,V)$ is an epi such that $r\circ e=\mathrm{id}_V$ and $e\circ r=p$. For every subobject $V\prec U$ we choose a representative retract $(V,e_{V\prec U},r_{V\prec U})$, and retract morphisms are represented graphically as in figure~\ref{fig:retract}.

We will need the notion of Picard group of a monoidal category.
An object $J$ of $\mathcal{C}$ is called invertible, or a simple current, if there exists an object $J'$ such that $J\otimes J'\cong\mathbf{1}$. One shows that an invertible object $J$ is simple, that $J\otimes J^\vee\cong\mathbf{1}\cong J^\vee\otimes J$, and that $\mathrm{dim}(J)=\pm1$. Since we assume the dimension of any object is positive, the dimension of any invertible object must be $1$.
The Picard group, or simple current group, of $\mathcal{C}$, $\mathrm{Pic}(\mathcal{C})$, is the multiplicative subgroup of the Grothendieck ring $K_0(\mathcal{C})$ generated by isomorphism classes of invertible objects. In particular we have an identification of a subset $G\subset I$ with $\mathrm{Pic}(\mathcal{C})$, and we use multiplicative notation so for $g,h\in G$, $gh\in G$ labels the isomorphism class of the object $U_g\otimes U_h$. Since $\mathcal{C}$ is braided, $\mathrm{Pic}(\mathcal{C})$ is abelian.
 Note that if $J$ is invertible and $U_i$ is any simple object, $J\otimes U_i$ is also simple and therefore isomorphic to $U_j$ for some $j\in I$. 
 Obviously $\mathrm{Pic}(\mathcal{C})$ acts on the isomorphism classes of simple objects, and hence on the set $I$. Write this action multiplicatively, so $U_g\otimes U_i\cong U_{gi}$ for $g\in G$, $i\in I$.
We also define the effective center $\mathrm{Pic}^\circ(\mathcal{C})\subset\mathrm{Pic}(\mathcal{C})$, a subset (not necessarily closed under multiplication) given by
\begin{equation}\eqlabel{efcenter}
\mathrm{Pic}^\circ(\mathcal{C})=\{g\in\mathrm{Pic}(\mathcal{C})|\  \theta_g^{|g|}=1\},
\end{equation}
where $|g|$ is the order of the element $g$.
In RCFT, the twist $\theta_i$ of a simple object $U_i$ is related to the conformal dimension $\Delta_i$ of the corresponding object (primary field) through
\begin{equation}\eqlabel{cdim}
\theta_i=\exp(-2\pi i \Delta_i).
\end{equation}
An invertible object $J$ thus defines an element of the effective center iff 
\begin{equation}\eqlabel{scineffcenter}
|J|\Delta_J\in\mathbb{Z}.
\end{equation}

\paragraph{Reshetikhin-Turaev TQFT and Quantum Representations} Using a modular category $\mathcal{C}$ it is possible to construct a three-dimensional TQFT, $\mathbf{tqft}_\mathcal{C}$. We briefly review some aspects of these TQFT's, and refer to~\cite{T} for further details. Essentially, $\mathbf{tqft}_\mathcal{C}$ is a monoidal functor from a category of three-dimensional extended cobordisms to the category $\mathcal{V}ec_f(\mathbb{C})$ of finite-dimensional $\mathbb{C}$-vector spaces. An object of the geometric category is an extended surface, an oriented closed topological surface $\Sigma$ with a finite (posibly empty) set of marked arcs, each labeled by an object of $\mathcal{C}$ and a sign, plus a choice of Lagrangian subspace $\lambda_\Sigma\subset H_1(\Sigma,\mathbb{R})$. We use the same symbol $\Sigma$ to refer to an extended surface as well as an underlying topological surface. A morphism is an extended cobordism, i.e. an oriented three-manifold with an embedded (possibly empty) oriented ribbon graph.
Strands of the ribbon graph are oriented, have a core orientation and are labeled by objects of $\mathcal{C}$ while coupons are oriented and 
labeled by morphisms of $\mathcal{C}$ in the obvious way consistent with the labeling of ribbons. Ribbons may pierce the boundary, and then gives rise to marked arcs in the obvious way. Extended cobordisms form a monoidal category with tensor unit $\emptyset$ and tensor product given by disjoint union. A TQFT is a ``projective'' monoidal functor $\mathbf{tqft}_\mathcal{C}$, by which we mean that it fails to be a monoidal functor only in the following way. If $\mathrm{M}_1:\Sigma_1\rightarrow \Sigma_2$ and $\mathrm{M}_2:\Sigma_2\rightarrow \Sigma_3$ are two extended cobordisms, then
$$\mathbf{tqft}_\mathcal{C}(\mathrm{M}_2)\circ\mathbf{tqft}_\mathcal{C}(\mathrm{M}_1)=\kappa_\mathcal{C}^m\mathbf{tqft}_\mathcal{C}(\mathrm{M}_2\circ\mathrm{M}_1),$$
where
$$\kappa_\mathcal{C}=\frac{\sum_{i\in I}\theta_i^{-1}\mathrm{dim}(U_i)^2}{\left(\sum_{i\in I}\mathrm{dim}(U_i)^2\right)^{1/2}},$$
$m$ is the Maslov index of three Lagrangian subspaces of $H_1(\Sigma_2,\mathbb{R})$ constructed from $\lambda_{\Sigma_1}$, $\lambda_{\Sigma_3}$, and $\lambda_{\Sigma_2}$.
See Section~IV.3--IV.4 of~\cite{T} for details on the Maslov index. It is not difficult to show that if $\mathrm{M}$ is an extended cobordism that as a $3$-manifold is just a cylinder over a surface, and $\mathrm{M}:\Sigma\rightarrow\Sigma$ as an endomorphism of an extended surface $\Sigma$, then $\mathbf{tqft}_\mathcal{C}(\mathrm{M})\circ\mathbf{tqft}_\mathcal{C}(\mathrm{N})=\mathbf{tqft}_\mathcal{C}(\mathrm{M}\circ \mathrm{N})$, and $\mathbf{tqft}_\mathcal{C}(\mathrm{N}')\circ\mathbf{tqft}_\mathcal{C}(\mathrm{M})=\mathbf{tqft}_\mathcal{C}(\mathrm{N}'\circ\mathrm{M})$ for any compatible $\mathrm{N}$, $\mathrm{N}'$. In other words, the relevant Maslov index vanish.
For a fixed category $\mathcal{C}$ we use the notation $\mathcal{H}(\Sigma)=\mathbf{tqft}_\mathcal{C}(\Sigma)$ for the vector space associated to an extended surface. Henceforth we restrict to extended surfaces with no marked points, and by further abuse of notation we then write $\mathcal{H}_g$ for $\mathcal{H}(\Sigma)$ if the genus of $\Sigma$ is $g$. Note that if we apply the TQFT functor to a cobordism from $\emptyset$ to a genus $g$ surface, the result is a linear map from $\mathbb{C}$ to $\mathcal{H}_g$. We occasionally apply such maps to $1\in\mathbb{C}$ to define elements in $\mathcal{H}_g$. The functor $\mathbf{tqft}_\mathcal{C}$ also has the property that it coincides on any two $\mathrm{M}$ and $\mathrm{M}'$ related by either an orientation preserving homeomorphism $f$ that restricts to the identity on $\partial\mathrm{M}=\partial\mathrm{M}'$ (in particular it only depends on the isotopy class of the ribbon graph), or a ``local move'' of the ribbon graph. Consider a part of an embedded ribbon graph that can be flattened, and hence be interpreted as a morphism in $\mathcal{C}$. A local move on that part of the ribbon graph is the procedure of replacing it with any other ribbon graph representing the same morphism in $\mathcal{C}$.

\begin{remark}
Let $\mathcal{C}$ be a modular category, and let $\Sigma$ be a closed oriented surface of genus $g$. Every oriented simple closed curve $\alpha$ in $\Sigma$ gives a linear representation $\rho^\alpha$ of $\mathrm{Pic}(\mathcal{C})$ on $\mathcal{H}_g$ defined by $\rho^\alpha(h)=\mathbf{tqft}_\mathcal{C}(\mathrm{M}_\Sigma^{\alpha, h})$, where $\mathrm{M}_\Sigma^{\alpha,h}$ is the cylinder over $\Sigma$ with a single ribbon labeled by $U_h$ running along $\alpha$, with no twist (zero framing) w.r.t. the boundary, in some horisontal section in the interior of the cylinder. If $\alpha$ is null homotopic, then $\rho^\alpha$ is the trivial representation. If $\alpha$ and $\beta$ are two oriented simple closed curves with geometric intersection number $0$, the maps $\rho^\alpha(h)$ and $\rho^\beta(h)$ commute, and in particular the elementwise composition defines a new representation $\rho^\alpha\cdot\rho^\beta$. Clearly $\rho^\alpha=\rho^{\alpha'}$ if $\alpha$ and $\alpha'$ are homotopic. Similarly any collection $S=\{\alpha_1,\ldots,\alpha_n\}$ of mutually disjoint oriented simple closed curves define a representation $\rho^S$ of $\mathrm{Pic}(\mathcal{C})$ on $\mathcal{H}_g$ such that $\rho^{\{\alpha,\beta\}}=\rho^\alpha\cdot\rho^\beta$ when $\alpha$ and $\beta$ are disjoint. When two oriented simple closed curves $\alpha$ and $\beta$ intersect only in a finite number of points we have the identity $\rho^\alpha(h)\circ\rho^\beta(h)=\theta_h^{\hat{\imath}(\alpha,\beta)}\rho^{\alpha\sharp\beta}(h)$ where $\hat{\imath}(\alpha,\beta)$ is the algebraic intersection number and $\alpha\sharp\beta$ is the collection of disjoint oriented simple closed curves obtained by smoothing each intersection of $\alpha$ and $\beta$ in the unique way consistent with orientations. This generalises a structure in the $\mathfrak{su}(2)$ case that was used in~\cite{BHMV} to construct an action of a Heisenberg group on $\mathcal{H}_g$,  which in turn is responsible for the reducibility of the corresponding quantum representations shown in~\cite{AM}.
\end{remark}

A functor $\mathbf{tqft}_\mathcal{C}$ gives a projective representation of the mapping class group of any surface underlying an extended surface via the mapping cylinder construction. Let $\Sigma$ be an extended surface with no marked points, and $f$ an orientation preserving homeomorphism of $\Sigma$. Then
\begin{equation}\eqlabel{mapcyl}
\mathrm{M}_f:=\Sigma\times[-1,0]\circ_f \Sigma\times[0,1]
\end{equation}
is the mapping cylinder of $f$, where $\circ_f$ means gluing by identification of points using $f$. The functor $\mathbf{tqft}_\mathcal{C}$ defines a linear map
\begin{equation}
\rho_g(f):=\mathbf{tqft}_\mathcal{C}(\mathrm{M}_f):\mathcal{H}_g\rightarrow\mathcal{H}_g,
\end{equation}
where the dependence on $\mathcal{C}$ is supressed.
The linear map $\rho_g(f)$ only depends on the homotopy class of $f$, and gives a projective representation of $MCG(\Sigma)$. We refer to the representations $V_g=(\mathcal{H}_g,\rho_g)$ as quantum representations.

The index set $I$ defines a basis $\{v_i\}_{i\in I}$ of $\mathcal{H}_1$ as follows. Let $\mathrm{H}_i$ be an oriented genus $1$ handlebody with a ribbon labeled by the object $U_i$ running along the non-contractible cycle without twist w.r.t. $\partial \mathrm{H}_1$. Then we can define $v_i:=\mathbf{tqft}_\mathcal{C}(\mathrm{H}_i)1\in\mathcal{H}_1$ (alternatively, define $v_i:=b_{U_i}\in\mathcal{H}_1$). 

\begin{figure}
\begin{picture}(155,90)
	\put(0,0)				{\includegraphics[scale=0.4]{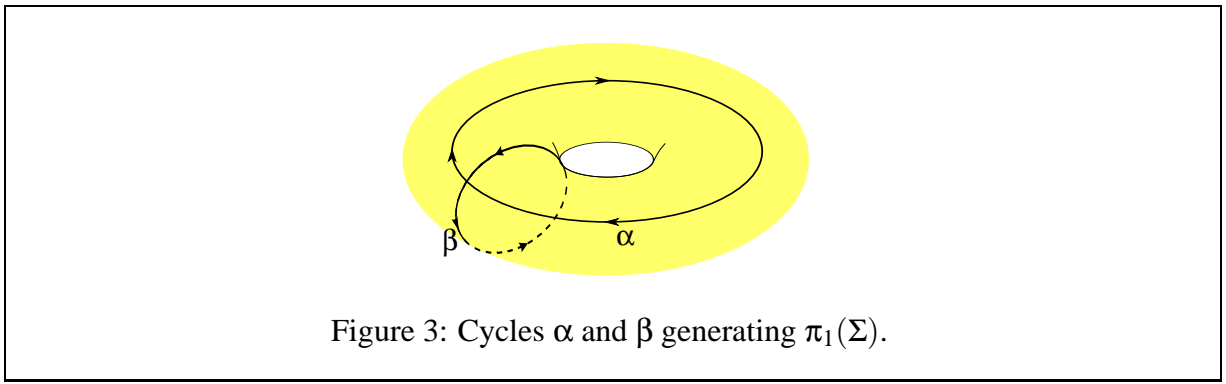}}
	\put(80,12)			{$\alpha$}
	\put(15,10)			{$\beta$}
\end{picture}
\caption{Cycles $\alpha$ and $\beta$ generating $\pi_1(\Sigma)$.}\label{fig:abcycles}
\end{figure}

Fix an extended surface $\Sigma$ of genus $1$, and choose oriented simple closed curves $\alpha$ and $\beta$ with geometric intersection number $1$ (see figure~\ref{fig:abcycles}).
The genus $1$ mapping class group is generated by the Dehn twists $T_\alpha$ and $T_\beta$ along $\alpha$ resp. $\beta$. The matrices representing $\rho_1(T_\alpha)$ respectively $\rho_1(T_\beta)$ in the basis $\{v_i\}_{i\in I}$ are easily determined by identifying $\bar v_i$ with $\mathbf{tqft}_{\mathcal{C}}(-\mathrm{H}_i)$, resulting in $\rho_1(T_\alpha)_{i,j}=\frac{s_{i,j}}{s_{0,0}}$, $\rho_1(T_\beta)_{i,j}=\theta_{i}\delta_{i,j}$.
One checks that the $|I|\times|I|$ matrices
$$S_{i,j}=\frac{s_{i,j}}{s_{0,0}},\quad T_{i,j}=\zeta_\mathcal{C}^{-1}\theta_i\delta_{i,j}$$
satisfy $(ST)^3=S^2$, $S^4=1$, 
where $$\zeta_\mathcal{C}=\left(\frac{\sum_{i\in I}\theta_i(\mathrm{dim}(U_i))^2}{\sum_{i\in I}\theta^{-1}_i(\mathrm{dim}(U_i))^2}\right)^{1/6}.$$
It furthermore turns out always to be the case that $S^2=C$, the charge conjugation matrix with elements $C_{i,j}=\delta_{\bar\imath,j}$.
Hence the genus $1$ representations lift to linear representations of $SL(2,\mathbb{Z})$, this is however not the case in general for higher genus.

 \section{A Reducibility Criterion}\selabel{redcrit}
We first review Theorem 1 of~\cite{AF}, that gives a condition for reducibility of quantum representations, and then show that the conditions of the theorem is satisfied for the modular categories of integrable highest weight representations of $\widehat{\mathfrak{su}}(N)$ for $N>2$ and any level $k\in \mathbb{N}$.
The statement of the theorem requires some preparation. We start by discussing algebras in monoidal categories, and continue to discuss some aspects of RCFT in terms of TQFT. See the series of papers~\cite{FRS1, FRS2, FRS3, FRS4, FjFRS1, FjFRS2, FrFRS1, FrFRS2} for a much more extensive discussion of algebras in monoidal categories and RCFT.
 
\begin{figure}
	\begin{picture}(400,155)(0,5)
		\put(30,100)				{\includegraphics[scale=0.3]{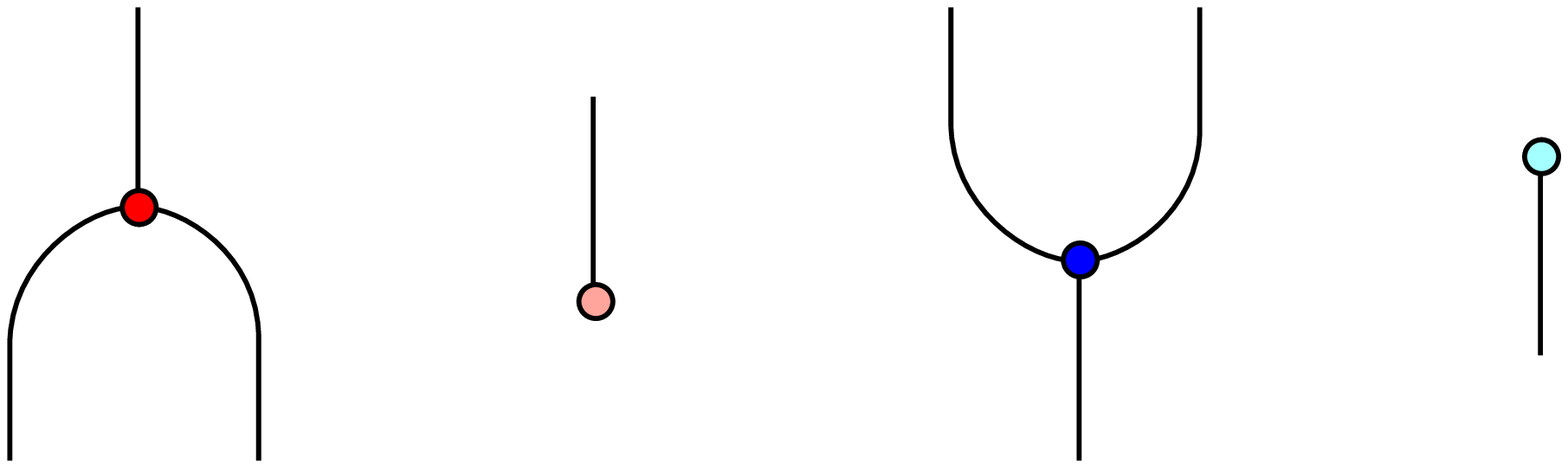}}
		\put(5,120)				{$m=$}
		\put(65,120)				{$\eta=$}
		\put(105,120)				{$\Delta=$}
		\put(160,120)				{$\varepsilon=$}
		\put(25,90)				{$A$}
		\put(53,90)				{$A$}
		\put(41,150)				{$A$}
		\put(86,140)				{$A$}
		\put(134,90)				{$A$}
		\put(123,150)				{$A$}
		\put(147,150)				{$A$}
		\put(180,100)				{$A$}
		\put(0,85)					{\line(1,0){400}}
		\put(210,0)				{\line(0,1){155}}
		\put(0,25)					{\includegraphics[scale=0.22]{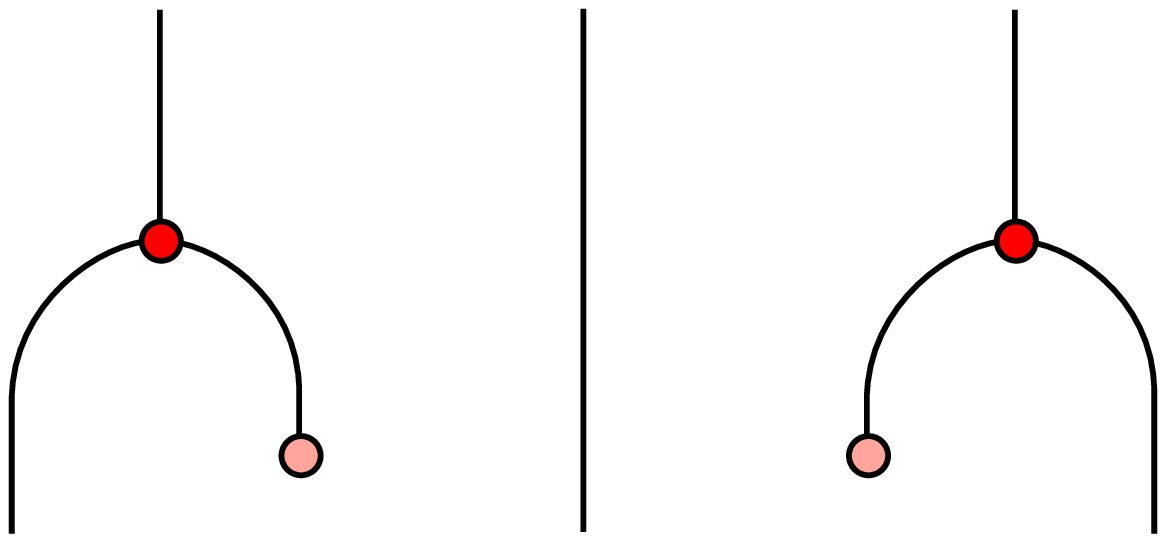}}
		\put(22,40)				{$=$}
		\put(42,40)				{$=$}
		\put(110,25)				{\includegraphics[scale=0.22]{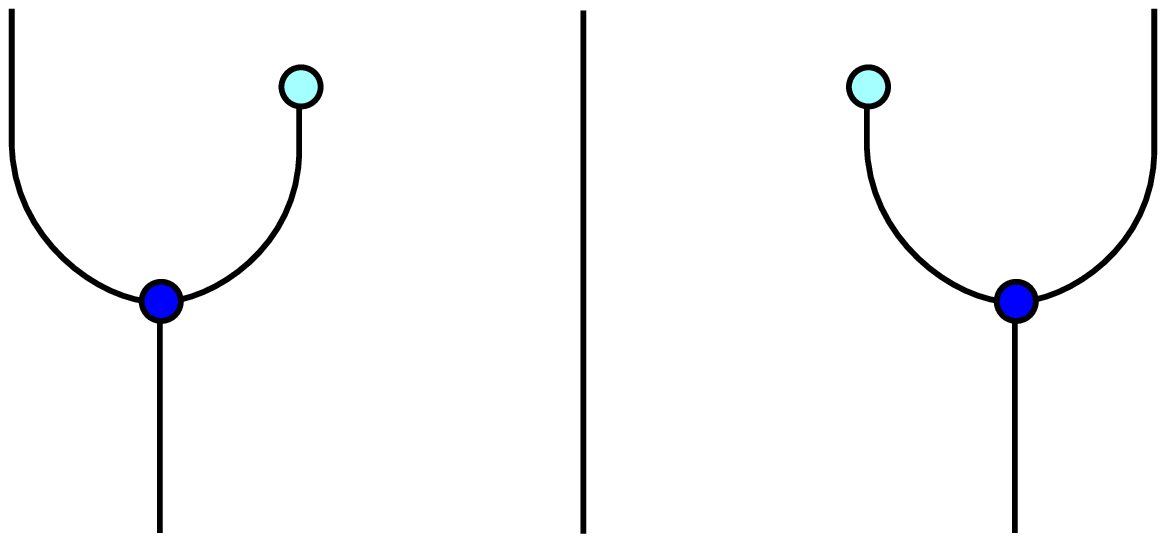}}
		\put(132,40)				{$=$}
		\put(152,40)				{$=$}
		\put(240,90)				{\includegraphics[scale=0.27]{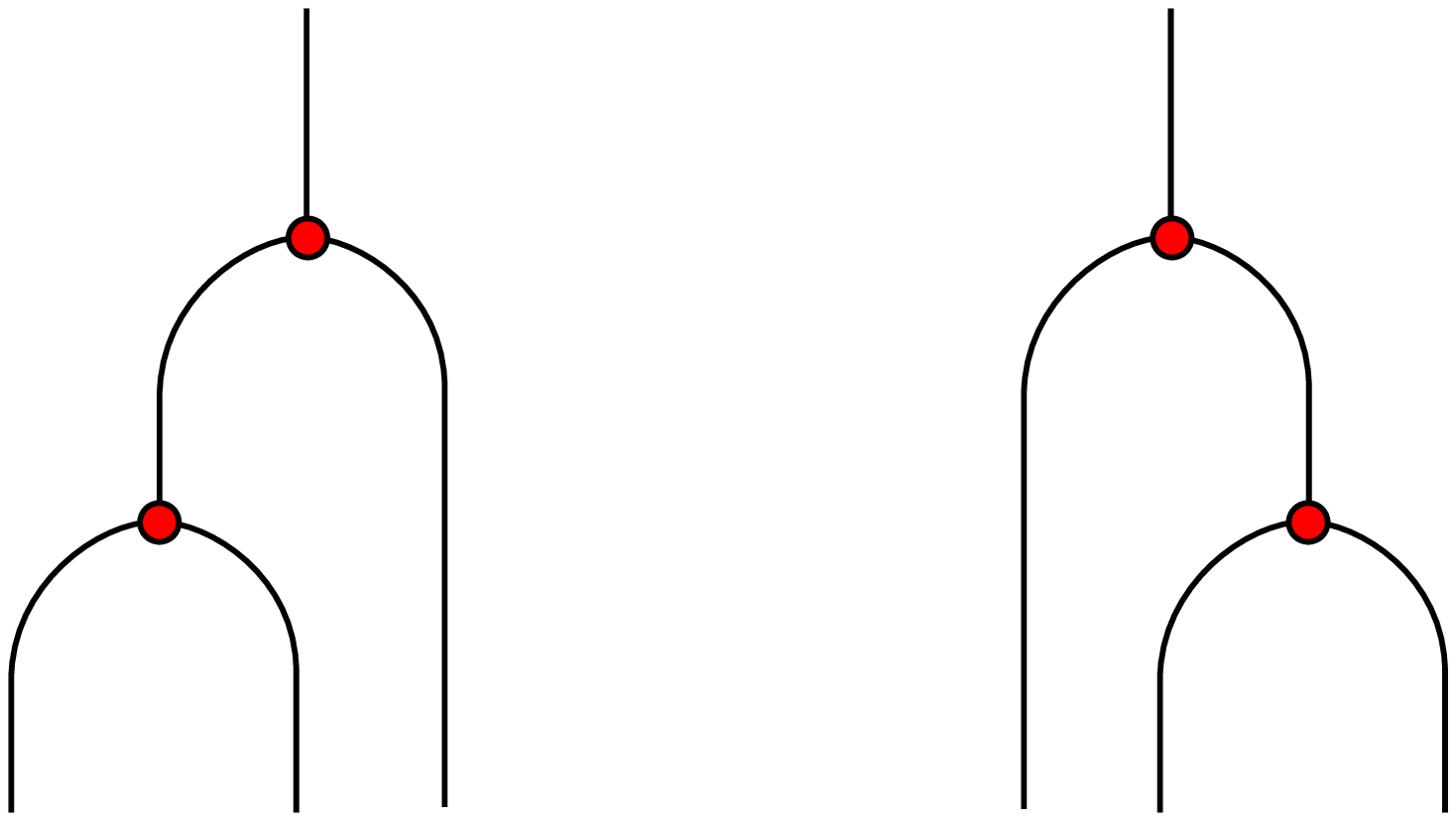}}
		\put(293,120)				{$=$}
		\put(240,5)				{\includegraphics[scale=0.27]{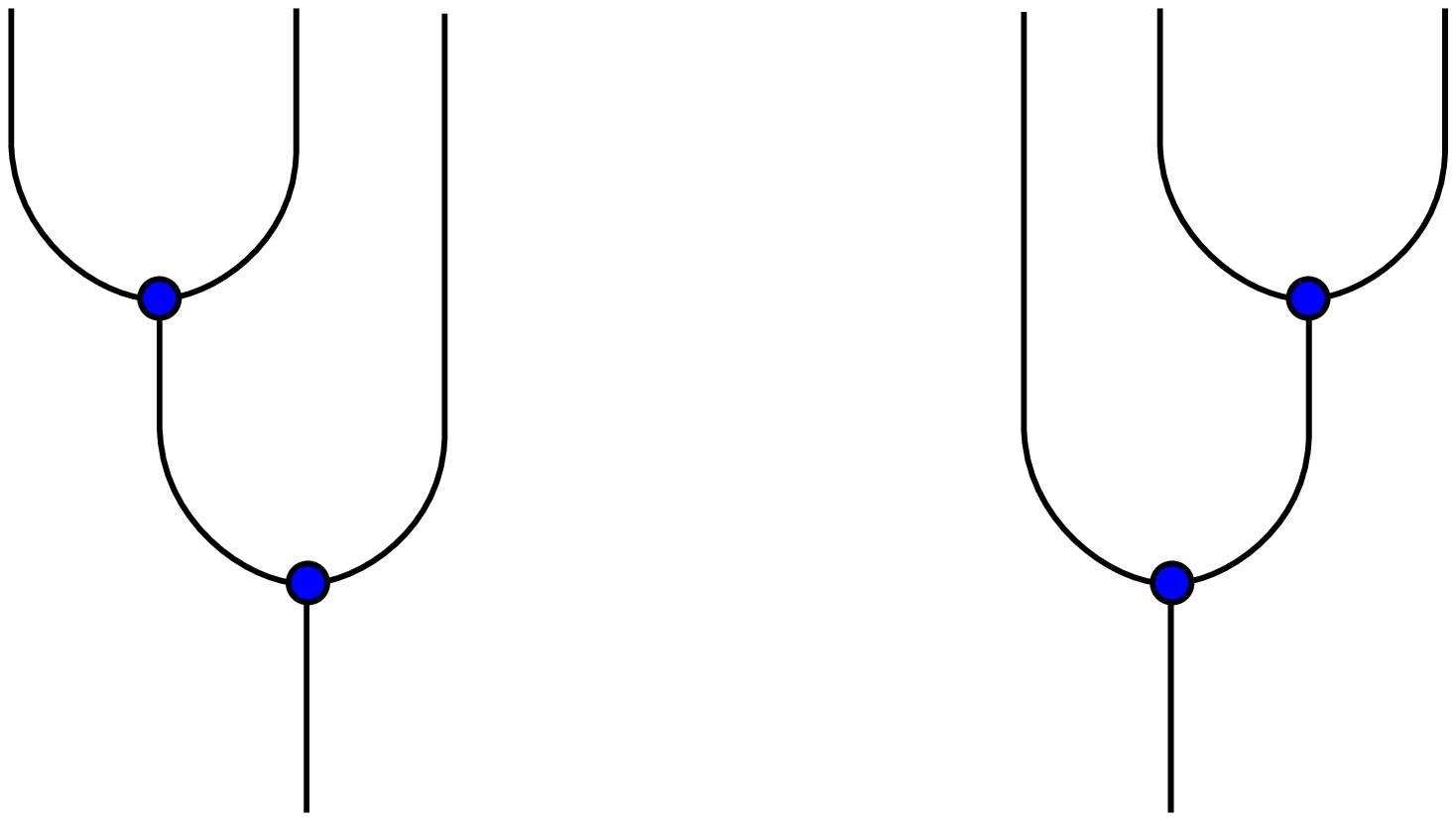}}
		\put(293,35)				{$=$}
	\end{picture}
\caption{Graphical representation of (co-)multiplication and (co-)unit morphisms plus (co-)unitality and (co-)associativity conditions for an algebra resp. coalgebra in $\mathcal{C}$.}
\label{fig:algmorph}
\end{figure}

\begin{figure}
	\begin{picture}(240,95)
		\put(0,0)				{\includegraphics[scale=0.4]{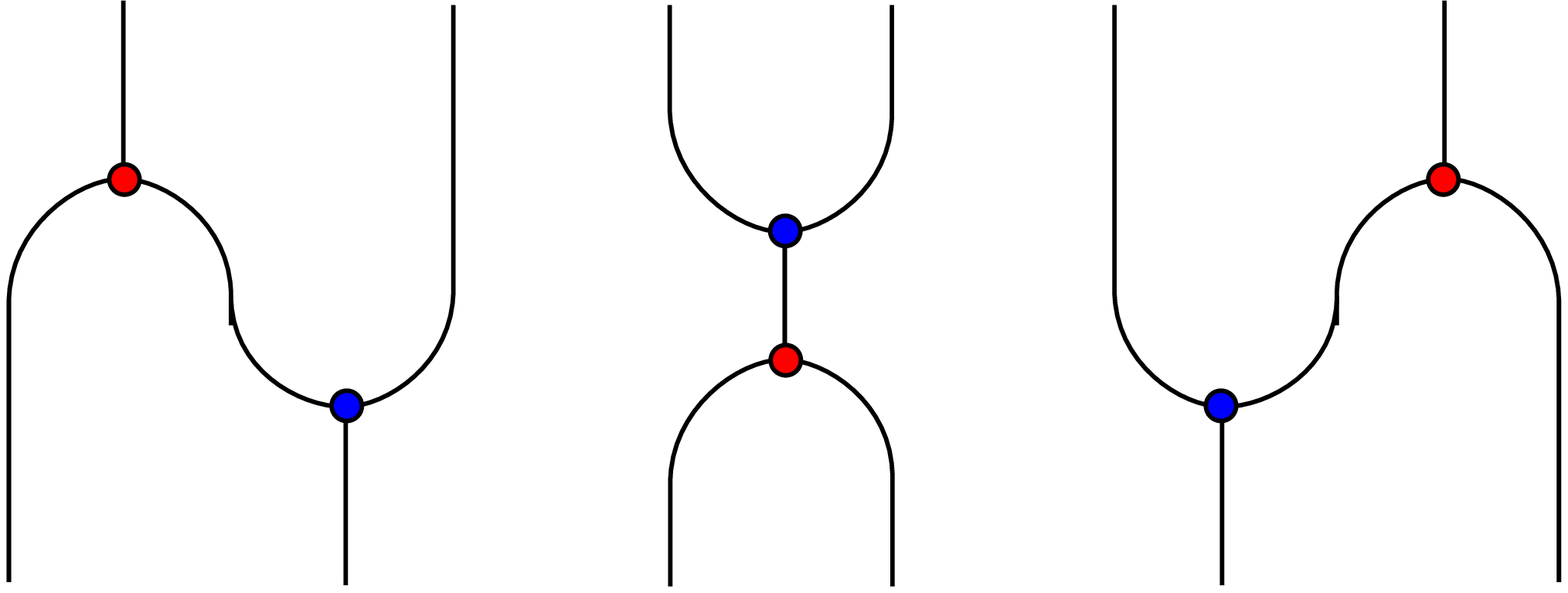}}
		\put(83,42)			{$=$}
		\put(140,42)			{$=$}
	\end{picture}
\caption{The Frobenius condition for a tuple $(A,m,\Delta)$.}
\label{fig:frob}
\end{figure}

\begin{figure}
\begin{picture}(400,90)
	\put(0,7)				{\includegraphics[scale=0.3]{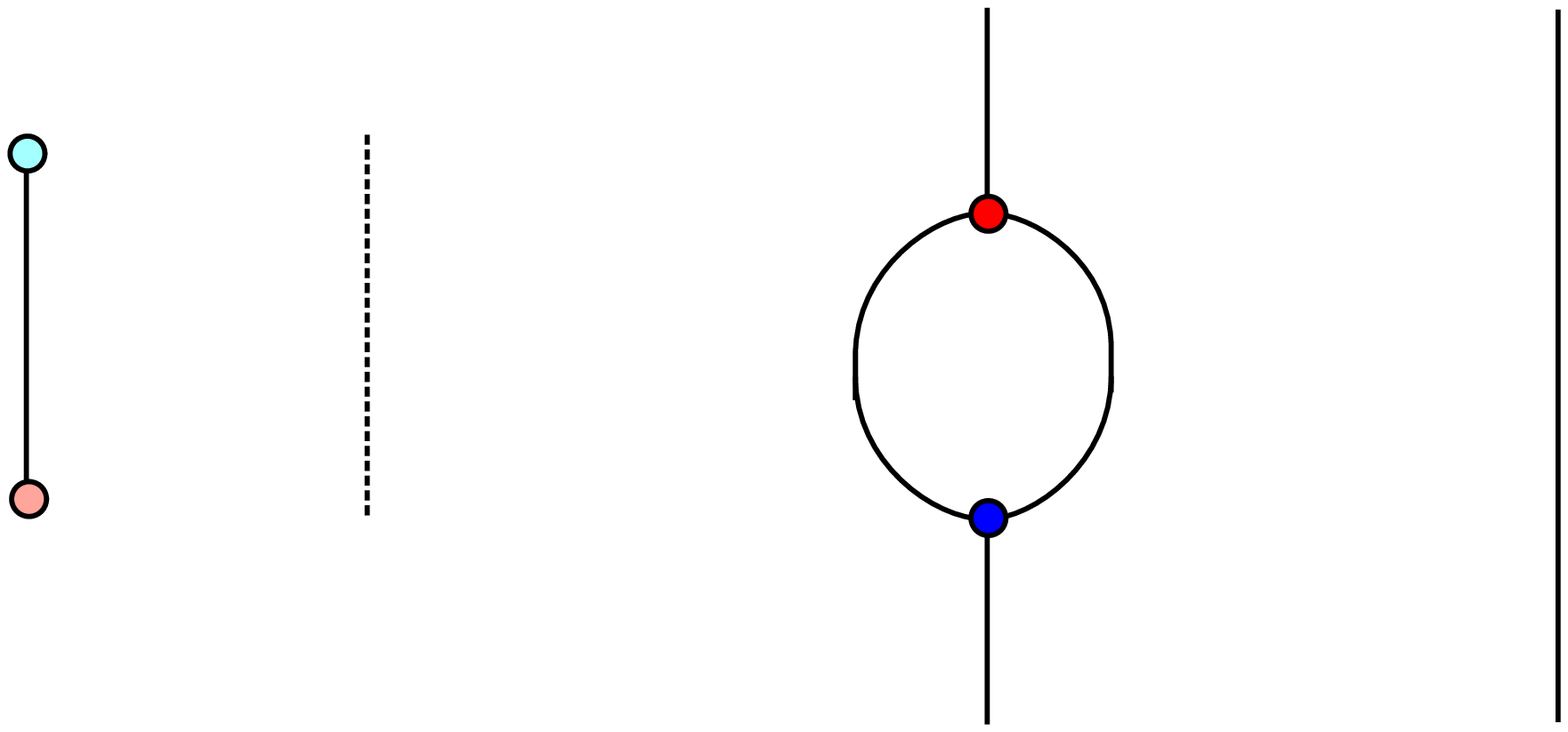}}
	\put(8,45)				{$=$\tiny$\ \mathrm{dim}(A)$}
	\put(122,42)			{$=$}
	\put(240,2)			{\includegraphics[scale=0.3]{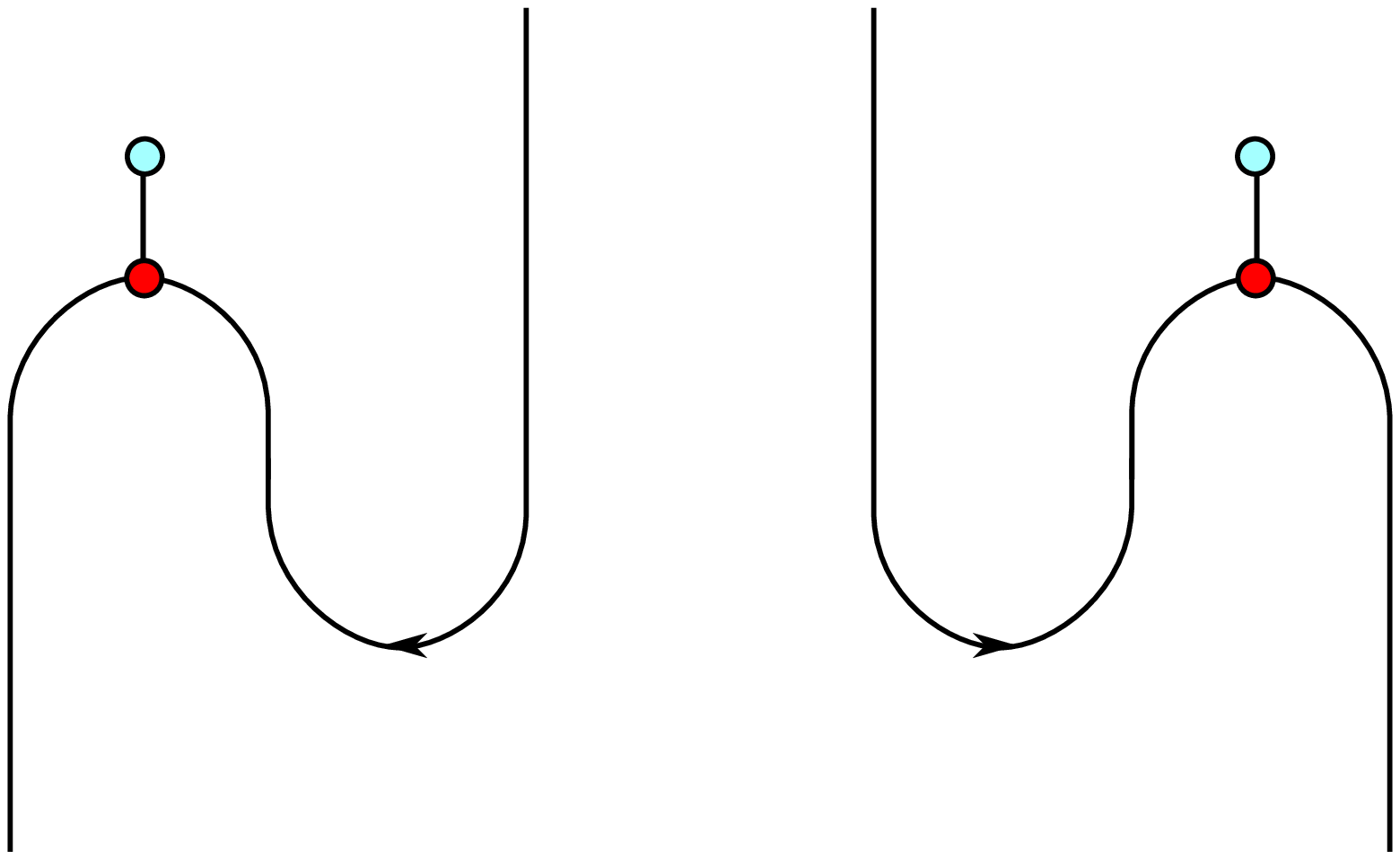}}
	\put(303,45)			{$=$}
\end{picture}
\caption{Graphical representation of the properties ``special'' and ``symmetric''.}\label{fig:specialsymm}
\end{figure}

 \paragraph{Algebras in Monoidal Categories}
An algebra in a monoidal category $\mathcal{C}$ is a triple $(A,m,\eta)$ where $A\in\mathrm{Obj}(\mathcal{C})$, $m\in\mathrm{Hom}(A\otimes A,A)$, $\eta\in\mathrm{Hom}(\mathbf{1},A)$ satisfying unitality $m\circ(\eta\otimes\mathrm{id}_A) = \mathrm{id}_A = m\circ(\mathrm{id}_A\otimes\eta)$ and associativity $m\circ(m\otimes\mathrm{id}_A)=m\circ(\mathrm{id}_A\otimes m)$. A coalgebra is defined analogously as a triple $(A,\Delta,\varepsilon)$, where the counit $\varepsilon\in \mathrm{Hom}(A,\mathbf{1})$ and the comultiplication $\Delta\in\mathrm{Hom}(A,A\otimes A)$ satisfy counitality resp. coassociativity conditions. (Co)multiplication and (co)unit morphisms are displayed in figure~\ref{fig:algmorph}, together with (co)unitality resp. (co)associativity conditions.
A Frobenius algebra is a tuple $(A,m,\eta,\Delta,\varepsilon)$ forming an algebra and a coalgebra, such that the Frobenius condition shown in figure~\ref{fig:frob} holds.
A Frobenius algebra $(A,m,\eta,\Delta,\varepsilon)$ is called special if $m\circ\Delta=\mathrm{id}_A$ and $\varepsilon\circ\eta=\mathrm{dim}(A)$\footnote{This is called normalised special in~\cite{FRS1, FRS2, FRS3, FRS4, FjFRS1, FjFRS2, FrFRS1, FrFRS2}.}, and it is called symmetric if the two obvious isomorphisms from $A$ to $A^\vee$ coincide. Graphical representations of specialness and symmetry are diplayed in figure~\ref{fig:specialsymm}.

One introduces (left) $A$-modules of an algebra $A$ in $\mathcal{C}$ as pairs $M=(\dot{M},\rho)$ where $\dot{M}\in\mathrm{Obj}(\mathcal{C})$ and $\rho\in\mathrm{Hom}(A\otimes\dot{M},\dot{M})$ satisfying the module properties $\rho\circ(\mathrm{id}_A\otimes\rho)=\rho\circ(m\otimes\mathrm{id}_{\dot{M}})$ and $\rho\circ(\eta\otimes\mathrm{id}_{\dot{M}})=\mathrm{id}_{\dot{M}}$. By defining module morphisms as morphisms in $\mathcal{C}$ intertwining the representation morphisms $\rho$ one gets the category $\mathcal{C}_A$ of left $A$-modules. Similarly one defines $A|B$-bimodules as triples $K=(\dot{K},\rho_L,\rho_R)$ where $\dot{K}\in\mathrm{Obj}(\mathcal{C})$, and $(\dot{K},\rho_L)$ is a left $A$-module while $(\dot{K},\rho_R)$ is a right $B$-module such that $\rho_L\circ(\mathrm{id}_A\otimes\rho_R)=\rho_R\circ(\rho_L\otimes\mathrm{id}_B)$. The category $\mathcal{C}_{A|B}$ of $A|B$-bimodules is defined by stating that bimodule morphisms are morphisms in $\mathcal{C}$ of the underlying objects that intertwine the left- and right-representation morphisms $\rho_L$ and $\rho_R$. An algebra $A$ in $\mathcal{C}$ is called simple if it is absolutely simple as a bimodule over itself, i.e. if $\mathrm{dim}_\mathbb{C}\mathrm{Hom}_{A|A}(A,A)=1$, and it is called haploid if $\mathrm{dim}_{\mathbb{C}}\mathrm{Hom}(\mathbf{1},A)=1$. A haploid algebra is automatically simple, and a haploid special algebra is in particular symmetric~\cite{FRS1}. Every monoidal category contains at least one symmetric special Frobenius algebra given by the object $\mathbf{1}$ where all morphisms are given by unit constraints.
Two algebras $A$ and $B$ are called Morita equivalent if $\mathcal{C}_A$ and $\mathcal{C}_B$ are equivalent (as module categories over $\mathcal{C}$~\cite{O}). Finally, it is known that any simple symmetric Frobenius algebra is Morita equivalent to a haploid symmetric Frobenius algebra~\cite{O}.

A class of algebras, central to the results of this paper, are the so-called Schellekens algebras.
A Schellekens algebra is a haploid special Frobenius algebra $A$ such that every simple subobject of $A$ is invertible. Schellekens algebras in modular categories were investigated in great detail in~\cite{FRS3}, with in particular the following results. Let $A$ be a Schellekens algebra in $\mathcal{C}$.
\begin{enumerate}
\item The support of $A$, $H(A):=\{g\in G|U_g\prec A\}$, is a group w.r.t. multiplication. (\cite{FRS3}, Lemma 3.13 (ii))
\item $H(A)\subset G$ corresponds to a subset of the effective center $\mathrm{Pic}^\circ(\mathcal{C})$ of $\mathcal{C}$.
\item Schellekens algebras in $\mathcal{C}$ are characterized by pairs $(H,\Xi)$ where $H$ is a subgroup of $\mathrm{Pic}^\circ(\mathcal{C})$ and $\Xi:H\times H\rightarrow\mathbb{C}^\times$ is a bihomomorphism such that $\Xi(h,h)=\theta_h$ for every $h\in H$, in the language of~\cite{FRS3} a Kreuzer-Schellekens bihomomorphism. To any such pair there exists a Schellekens algebra $A$ such that $H(A)\cong H$. (\cite{FRS3}, Corrollary 2.18, Remark 2.21 (i), Proposition 3.14 (ii))
\item On every subgroup of the effective center of $\mathcal{C}$ there exists at least one Kreuzer-Schellekens bihomomorphism, so there is at least one Schellekens algebra associated to every subgroup of the effective center. (\cite{FRS3}, Proposition 3.20 (iii))
\end{enumerate}

\paragraph{Elements of RCFT}
By a consistent orientation of a trivalent graph we mean an orientation of each edge such that every vertex is adjacent to either two incoming and one outgoing, or one incoming and two outgoing edges.
 Let $\Sigma$ be a closed oriented surface of genus $g$, let $T$ be a consistently oriented dual triangulation\footnote{By ``dual triangulation'' we mean an embedded graph whose Poincar\'e dual graph is a triangulation, in particular this implies that all vertices of $T$ are trivalent, and for every homotopy class of closed curves in $\Sigma$ there exist a cycle of $T$ representing that class. It is not difficult to show that every dual triangulation of a surface can be given a consistent orientation.} of $\Sigma$, and let $\lambda\subset H_1(\Sigma,\mathbb{R})$ be a Lagrangian subspace. Furthermore, let $A$ be a symmetric special Frobenius algebra in a modular category $\mathcal{C}$. Consider the oriented three-manifold $\Sigma\times[-1,1]$, with the orientation inherited from the orientations of $\Sigma$ and $[-1,1]$. Define an extended cobordism $\mathrm{M}[\Sigma,A,T]:\Sigma\rightarrow\Sigma$ by embedding a ribbon graph in $\Sigma\times[-1,1]$ as follows. In a neighbourhood of $T\subset\Sigma\times\{0\}\subset\Sigma\times[-1,1]$, cover edges of $T$ by planar ribbons and vertices by planar coupons. Provide the core of the ribbons with the orientation from $T$, and provide the faces of ribbons and coupons the opposite of the orientation of $\Sigma$. Label the ribbons with the object $A$, and the coupons by the multiplication (comultiplication) morphism of $A$ if the coupon covers a vertex with two incoming (outgoing) edges.
 \begin{figure}
 \begin{picture}(190,40)
 				\put(0,8)				{\includegraphics[scale=0.3]{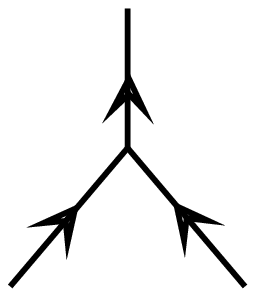}}
				\put(40,0)				{\includegraphics[scale=0.3]{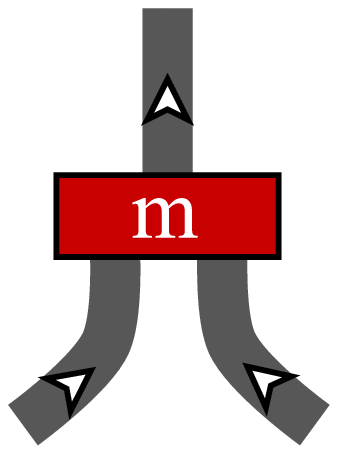}}
				\put(25,18)			{$\rightsquigarrow$}
				\put(58,30)			{\tiny $A$}
				\put(40,8)				{\tiny $A$}
				\put(67,8)				{\tiny $A$}
				\put(120,8)			{\includegraphics[scale=0.3]{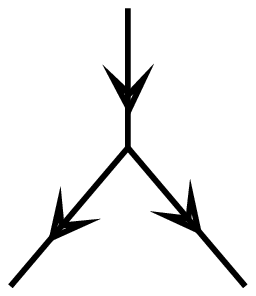}}
				\put(160,0)			{\includegraphics[scale=0.3]{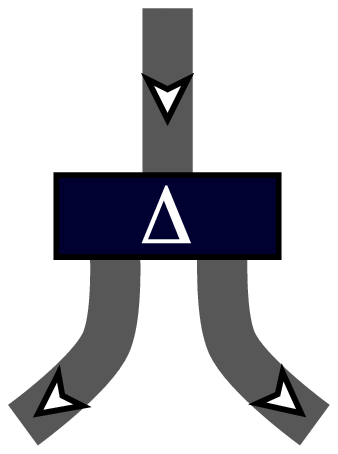}}
				\put(145,18)			{$\rightsquigarrow$}
				\put(178,30)			{\tiny $A$}
				\put(160,8)			{\tiny $A$}
				\put(187,8)			{\tiny $A$}
 \end{picture}
 \caption{Cover every vertex of $T$ by a coupon labeled by multiplication or comultiplication of the algebra $A$, depending on the number of incoming and outgoing edges.}
 \label{fig:vertexdecoration}
 \end{figure}

As an extended cobordism, $\mathrm{M}[\Sigma,A,T]$ can be hit with the functor $\mathbf{tqft}_\mathcal{C}$ to give a linear endomorphism of $\mathcal{H}(\Sigma)$.
\begin{definition}\delabel{correlator} Given $\mathrm{M}[\Sigma,A,T]$ as above, define the linear map 
$$\mathrm{P}[\Sigma,A,T]:=\mathbf{tqft}_\mathcal{C}(\mathrm{M}[\Sigma,A,T])\in\mathrm{End}_{\mathbb{C}}(\mathcal{H}(\Sigma)).$$
If $\Sigma$ has genus $g$ we abuse notation and refer to $\mathrm{P}[\Sigma,A,T]$ as the genus $g$ partition function associated to $A$ and $T$.
\end{definition}
It can be shown (Proposition 3.2 of~\cite{FjFRS1}) that $\mathrm{P}[\Sigma,A,T]$ is independent of the choice of dual triangulation $T$, we therefore drop reference to $T$ and write $\mathrm{P}[\Sigma,A]$, or only $\mathrm{P}_g[A]$, for the genus $g$ partition function associated to $A$. The proof essentially follows from the well-known fact that any two triangulations with the same number of faces are connected by a sequence of Whitehead moves, and for a consistently oriented dual triangulation a Whitehead move is a local move corresponding to either a (co-)associativity condition, or a Frobenis condition. Any two consistent orientations can be connected by local moves corresponding to the symmetry property, and the property ``special'' is used to add or remove edges of a dual triangulation.
The partition functions also have the property~\cite{FrFRS2} that if $A$ and $B$ are Morita equivalent algebras, then $\mathrm{P}_g[A]=\gamma^{g-1}\mathrm{P}_g[B]$ where $\gamma=\mathrm{dim}(A)/\mathrm{dim}(B)$. In particular, the genus $1$ partition functions of Morita equivalent algebras coincide.

\begin{figure}
\begin{picture}(100,110)
	\put(0,50)				{$P_A^l(U)=$}
	\put(50,10)			{\includegraphics[scale=0.4]{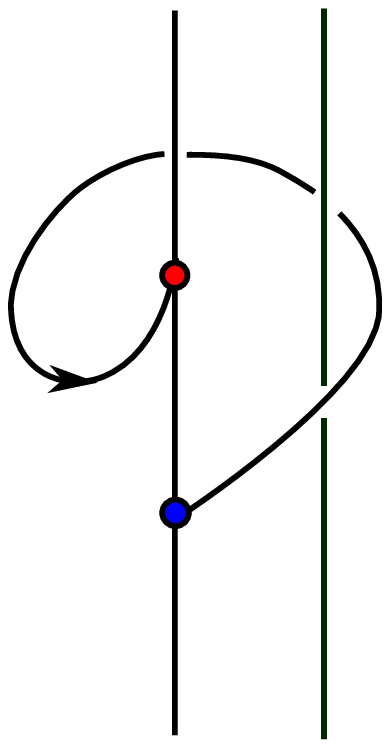}}
	\put(65,0)				{$A$}
	\put(82,0)				{$U$}
	\put(65,100)			{$A$}
	\put(82,100)			{$U$}
\end{picture}
\caption{The morphism $P_A^l(U)$ defining the object $E_A^l(U)$.}\label{fig:PlU}
\end{figure}

Given a special symmetric Frobenius algebra $A$ in a ribbon category $\mathcal{C}$ one defines an endofunctor $E^l_A$ of $\mathcal{C}$~\cite{FrFRS1} that acts on objects as
\begin{equation}\eqlabel{Efunct1}
U\in\mathrm{Obj}(\mathcal{C})\mapsto E^l_A(U)=\mathrm{Im}(P_A^l(U))
\end{equation}
where $P_A^l(U)\in\mathrm{End}(A\otimes U)$ is an idempotent defined in figure~\ref{fig:PlU}.
Using the properties symmetric, special and Frobenius, it is trivial to check that $P_A^l(U)$ is indeed an idempotent. The object $C_l(A):= E_A^l(\mathbf{1})$ carries the structure of a commutative symmetric special Frobenius algebra with structure morphisms defined by those of $A$ and the retract morphisms of $C_l(A)\prec A$, i.e. $m_{C_l(A)}=r_{C_l(A)\prec A}\circ m_A\circ (e_{C_l(A)\prec A}\otimes c_{C_l(A)\prec A})$ etc. This algebra is called the left center of $A$.
If the genus of $\Sigma$ is $1$ it is conventional to call the torus partition function $Z$, so we define $Z(A):=\mathrm{P}_1[A]$. If in the canonical basis $\{v_i\}_{i\in I}$ of $\mathcal{H}_1$ we write $Z(A)=\sum_{i,j\in I}Z_{i,j}(A)v_i\otimes\bar v_j$, where $\{\bar v_i\}_{i\in I}$ is the dual basis, one can show~\cite{FRS1} that the coefficients $Z_{i,j}(A)$ take values
\begin{equation}\eqlabel{toruspfn}
Z_{i,j}(A)=\mathrm{dim}_\mathbb{C}\mathrm{Hom}(E_A^l(U_i),U_j).
\end{equation}
For $A=\mathbf{1}$, the matrix is obviously $Z_{i,j}(\mathbf{1})=\delta_{i,j}$, and the linear map is the identity. It is equally simple to show that $P[\Sigma,\mathbf{1}]=\mathrm{id}_{\mathcal{H}(\Sigma)}$ for any $\Sigma$.

A torus partition function $Z(A)$ that is proportional to $\mathrm{id}_{\mathcal{H}_1}$ will be called {\em trivial}, and a non-trivial torus partition function is thus any $Z(A)$ not proportional to the identity.

We note that Schur's lemma implies reducibility of $V_1$ if there exists an $|I|\times |I|$-matrix $Z$ not proportional to the unit matrix and such that $SZ=ZS$, $TZ=ZT$. There is a sizeable amount of CFT literature devoted to find such matrices (satisfying some additional constraints such as $Z_{0,0}=1$, and $Z_{i,j}\in\mathbb{N}_0$ for any $i,j\in I$), and it is therefore possible to determine reducibility of certain genus one quantum representations. It is worth pointing out that if the category $\mathcal{C}$ is such that not all simple objects are self-dual, then the genus $1$ representation is always reducible since the charge conjugation matrix $C$ gives a non-trivial intertwiner.

For $g\in G$, $i\in I$ define 
\begin{equation}\eqlabel{char}
\chi_{U_i}(U_g)=\theta_{gi}\theta_g^{-1}\theta_i^{-1}.
\end{equation}
One can show~\cite{FRS3} that for each $i\in I$, the function $\chi_{U_i}:G\rightarrow\mathbb{C}^\times$ is a character of $G$.
A Schellekens algebra on the object $A=\oplus_{h\in H}U_h$ and with Kreuzer-Schellekens bihomomorphism $\Xi_A$ gives the following torus partition function~(Theorem 3.27 of \cite{FRS3})
\begin{equation}\eqlabel{scpnfn}
Z_{i,j}(A)=\frac{1}{|H|}\sum_{h,g\in H}\chi_{U_i}(U_h)\Xi_A(U_h,U_g)\delta_{j,gi}.
\end{equation}
The existence of a Schellekens algebra for every subgroup of the effective center of $\mathcal{C}$ together with the formula~\equref{scpnfn} will be crucial for showing the reducibility of quantum representations in the case of $\mathcal{C}_{N,k}$.

\begin{remark}
\begin{enumerate}
\item[(i)] The matrix $Z_{i,j}(A)$ is known as the modular invariant torus partition function in CFT because it commutes with the matrices $S$ and $T$ constructed from $\mathcal{C}$ that generate $SL(2,\mathbb{Z})$. Actually, in parts of the literature it is not this matrix, but the matrix $\widetilde{Z}(A):=Z(A)C$ that is called the torus partition function. In this convention, the trivial algebra $\mathbf{1}$ corresponds to the charge conjugation matrix $C$.
\item[(ii)] The torus partition functions \equref{scpnfn} are precisely the so-called simple current invariants from~\cite{KrS} multiplied by the matrix $C$.
\end{enumerate}
\end{remark}

\paragraph{Main Theorem}

The main result, first proven in~\cite{AF}, is the following theorem.
\begin{theorem}\thlabel{red}(Theorem 1 of \cite{AF})
 Let $\mathcal{C}$ be a modular category. If there exists a symmetric special Frobenius algebra $A$ in $\mathcal{C}$ such that the matrix $Z(A)$ is not proportional to the unit matrix, then the representations $V_g$ provided by $\mathbf{tqft}_\mathcal{C}$ are reducible for all $g\geq 1$.
 \end{theorem}
 
 We will give an outline of the proof, referring to~\cite{AF} for the remaining details. The idea behind the proof is very simple; provided the conditions of the theorem hold we show that the genus $g$ partition function $\mathrm{P}_g[A]$ is an intertwiner of $V_g$, not proportional to the identity, for every genus $g\geq 1$. The proof of the theorem then follows from Schur's lemma. Note that the fact that quantum representations are only projective, and not linear, representations is irrelevant for this argument. For the rest of this section, fix a modular category $\mathcal{C}$ and a symmetric special Frobenius algebra $A$ in $\mathcal{C}$. 

Let us first show that $\mathrm{P}_g[A]$ is a self-intertwiner of $V_g$.
 \begin{lemma}\lelabel{minv}
 Let $\Sigma$ be an extended surface of genus $g$ with no marked points, and let $f$ be an orientation preserving homeomorphism of $\Sigma$, then
 \begin{equation}\eqlabel{modinv}
 P[\Sigma,A]\circ\rho_{g}(f)=\rho_{g}(f)\circ P[\Sigma, A]
 \end{equation}
 \end{lemma}
 
The lemma follows essentially from Theorem 2.2 of~\cite{FjFRS1}, but here we provide a slightly different proof.
\begin{proof}
We note first that $\mathrm{P}[\Sigma,A]\circ\rho_g(f)=\kappa_\mathcal{C}^m \mathbf{tqft}_\mathcal{C}(\mathrm{M}[\Sigma,A,T] \circ\mathrm{M}_f)$, where $m$ is the relevant Maslov index. In fact, since $\mathrm{M}[\Sigma,A,T]$ as a three-manifold is just $\Sigma\times[-1,1]$, we know that $m=0$. Consider pushing $f$ ``down'' through the ribbon graph inside $\mathrm{M}[\Sigma,A,T]\circ\mathrm{M}_f$. The only thing that can happen as a result is that the triangulation $T$ changes to a different triangulation $T'$ so we end up with $\mathrm{M}_f\circ\mathrm{M}[\Sigma,A,T']$. Use again $\rho_g(f)\circ\mathrm{P}[\Sigma,A]=\kappa_\mathcal{C}^{m'}\mathbf{tqft}_\mathcal{C}(\mathrm{M}_f\circ\mathrm{M}[\Sigma,A,T'])$ (using the fact that the partition function is independent of $T$), and the same argument as above implies $m'=0$.
\end{proof}
 
 We have thus shown that every symmetric special Frobenius algebra gives intertwiners of the representations $\rho_g$, $g\in\mathbb{N}$. To prove theorem~\thref{red} we must show that $\mathrm{P}_g[A]$ is non-trivial for every $g\geq 1$ if $Z(A)$ is non-trivial.
 
\begin{proof} (of \thref{red}) We give a broad outline of the proof, and again refer to~\cite{AF} for details.
For an object $U$, let $v^g_U\in\mathcal{H}_g$ be the element given by acting with the functor $\mathbf{tqft}_\mathcal{C}$ on the genus $g$ extended handlebody $\mathrm{H}^g_U$ shown in figure~\ref{fig:vgi} and applying the resut to $1\in\mathbb{C}$.
\begin{figure}
\begin{picture}(300,70)
	\put(30,0)				{\includegraphics[scale=0.5]{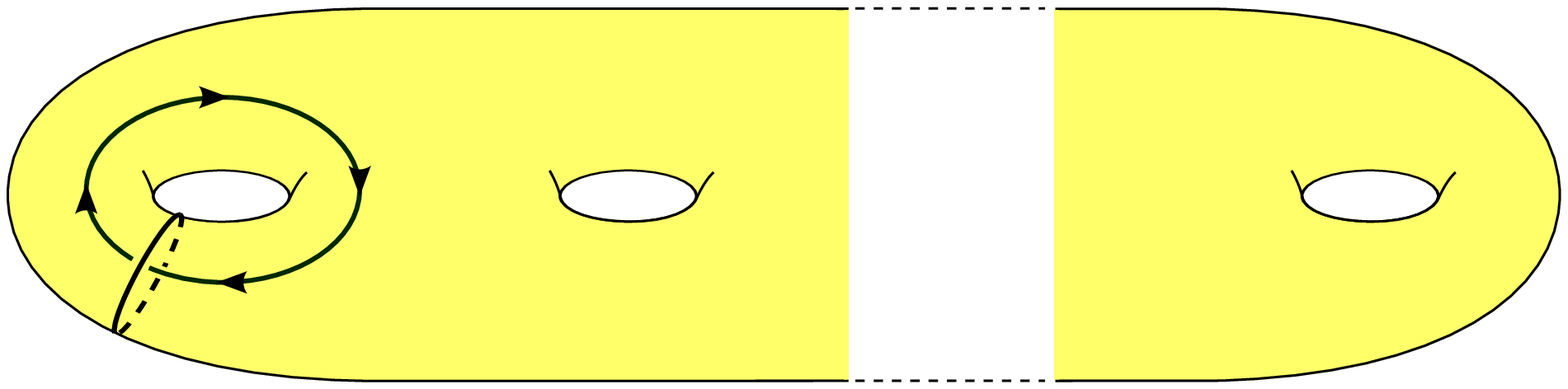}}
	\put(0,32)				{$\mathrm{H}^g_U=$}
	\put(48,50)			{$U$}
\end{picture}
\caption{The genus $g$ extended handlebody $\mathrm{H}^g_U$ with a single ribbon labeled by $U_i$. The ribbon has no twist (zero framing) w.r.t the boundary $\partial\mathrm{H}^g_U$.}\label{fig:vgi}
\end{figure}
When $U=U_i$ for some $i\in I$ we write $v_i^g\equiv v_{U_i}^g$ and $\mathrm{H}^g_i\equiv\mathrm{H}^g_{U_i}$.
In terms of the definition of $\mathcal{H}_g$ as a Hom-space of $\mathcal{C}$ (see~\cite{T}), the element $v^g_i$ is given by the tensor product of a number of unit constraints together with one instance of $b_{U_i}$, in particular $v^g_i\neq 0$. Choose the consistently oriented dual triangulation of the genus $g$ surface $\Sigma=\partial\mathrm{H}^g_i$ shown in figure~\ref{fig:dualtri} to construct $\mathrm{M}[\Sigma,A,T]$.
\begin{figure}
\begin{picture}(270,70)
	\put(0,0)				{\includegraphics[scale=0.5]{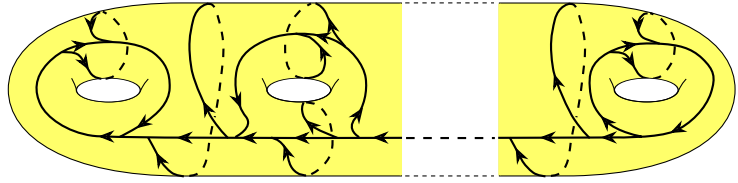}}
\end{picture}
\caption{Choice of consistently oriented dual triangulation of $\Sigma\cong\partial\mathrm{H}^g_i$.}\label{fig:dualtri}
\end{figure}

Up to a non-zero multiplicative factor, the result of applying $\mathrm{P}_g[A]$ to $v^g_i$ is obtained by applying $\mathbf{tqft}_\mathcal{C}$ to the extended cobordism obtained by gluing $\mathrm{M}[\Sigma,A,T]$ on top of $\mathrm{H}^g_i$, see figure~\ref{fig:PonH}.
\begin{figure}
\begin{picture}(285,70)
	\put(0,0)				{\includegraphics[scale=0.5]{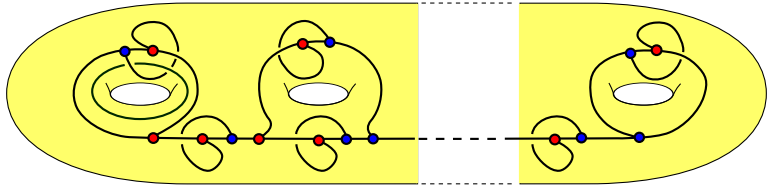}}
\end{picture}
\caption{The result of gluing $\mathrm{M}[\Sigma,A,T]$ onto $\mathrm{H}^g_i$. All ribbons except for one are labeled by $A$, and all trivalent morphisms are labeled by either $\Delta$ or $m$. The last ribbon is labeled by $U_i$.}\label{fig:PonH}
\end{figure}

Replace each occurence of $P_A^l(U)$ inside the resulting three-manifold with $e_{E^l_A(U)\prec A\otimes U}\circ r_{E^l_A(U)\prec A\otimes U}$, and slide the retract morphisms towards the coupons labeled by $m$ and $\Delta$. The result is a ribbon graph labeled by the left center $C_l(A)$, and by the object $E_A^l(U_i)$ together with the morphism shown in figure~\ref{fig:PonHmorph}.
\begin{figure}
\begin{picture}(70,120)
	\put(0,10)				{\includegraphics[scale=0.4]{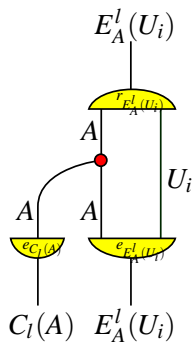}}
	\put(0,0)				{$C_l(A)$}
	\put(32,0)				{$E^l_A(U_i)$}
	\put(32,113)			{$E^l_A(U_i)$}
	\put(5,32)				{\tiny $e_{C_l(A)}$}
	\put(40,32)			{\tiny $e_{E^l_A(U_i)}$}
	\put(40,88)			{\tiny $r_{E^l_A(U_i)}$}
	\put(3,43)				{$A$}
	\put(27,43)			{$A$}
	\put(59,55)			{$U_i$}
	\put(27,73)			{$A$}
\end{picture}
\caption{Apart from multiplications and comultiplications of the left center $C_l(A)$, the remaining morphism in the extended cobordism resulting from gluing $\mathrm{M}[\Sigma,A,T]$ to $\mathrm{H}^g_i$, replacing all projectors by retract morphisms and sliding these toward the trivalent morphisms.}\label{fig:PonHmorph}
\end{figure}

A decomposition of $C_l(A)$ in simple objects correspond to a decomposition of $\mathrm{P}_g[A](v^g_i)$. We can without loss of generality assume $A$ is simple (every algebra is equivalent to a direct sum of simple algebras~\cite{O}), and even a haploid one (choosing a haploid representative of the original simple algebra). The last step of the proof amounts to showing that restricting to the subobject $\mathbf{1}\prec C_l(A)$ gives $\lambda v_{E^l_A(U_i)}^g$ for some $\lambda\in\mathbb{C}^\times$ independent of $i$.  
This is a straightforward exercise in calculus in a semisimple ribbon category (see~\cite{AF} for details).
By the definition of $\mathcal{H}_g$ we conclude that $\mathrm{P}_g[A](v_i^g)=\lambda v_{E^l_A(U_i)}^g+R_i^g$, where $R_i^g$ is either $0$ or  linearly independent of $\{v_i^g\}_{i\in I}$.  Using \equref{toruspfn} this implies $\mathrm{P}_g[A]\not\propto\mathrm{id}_{\mathcal{H}_g}$ if $Z(A)\not\propto\mathrm{id}_{\mathcal{H}_1}$.
\end{proof}
 
\section{The $SU(N)$ Case}\selabel{sun}
Let $\mathcal{C}_{N,k}$ be the category of integrable highest weight representations of $\widehat{su}(N)$ at level $k\in\mathbb{N}$, which is known to be modular~\cite{F, HL, BK}. Alternatively we can take $\mathcal{C}_{N,k}$ to be the modular category obtained by taking the semisimple subquotient of the category of representations of the quantum group $U_q(\mathfrak{su}(N))$ for $q$ a  $k+N$'th root of unity. We let $\lambda=(\lambda_1,\ldots,\lambda_{N-1})$ be the dynkin labels of a representation of $\mathfrak{su}(N)$, and define 
$$\mathcal{P}^{N,k}_{++}:=\{\lambda=(\lambda_1,\ldots,\lambda_{N-1})|\sum_{i=1}^{N-1}\lambda_i\leq k\}.$$
The elements of $\mathcal{P}_{++}^{N,k}$ label (isomorphism classes of) simple objects of $\mathcal{C}_{N,k}$.
We will use the existence of Schellekens algebras as described in \seref{redcrit} together with equation \equref{scpnfn} to assert the existence of symmetric special Frobenius algebras with non-trivial torus partition function in the category $\mathcal{C}$ for various $N$ and $k$.
The data we need to construct such an algebra is (i) a subgroup $H$ of the effective center of $\mathcal{C}$, and (ii) a Kreuzer-Schellekens bihomomorphism on $H$.
The Picard group $\mathrm{Pic}(\mathcal{C})$ is $\mathbb{Z}/N\mathbb{Z}$~\cite{SY}, and is generated by the simple object
\begin{equation}
	J=  (k,0,0,\ldots,0).
\end{equation}
We have that $J$ acts on simple objects as 
\begin{equation}\eqlabel{scaction}
J:(\lambda_1,\ldots,\lambda_{N-1})\mapsto(k-\Sigma_{i=1}^{N-1}\lambda_i,\lambda_1,\lambda_2,\ldots,\lambda_{N-2}),
\end{equation}
implying
$$J^p=(0,\ldots,0,k,0,\ldots,0)$$ with the only non-zero entry in the $p$'th position.
According to \equref{cdim}, to determine $\mathrm{Pic}^\circ(\mathcal{C})$ it is enough to determine the conformal weights of the simple objects $J^p$, $p=1,\ldots,N-1$. The conformal weight of a simple object $\lambda$ is given by
\begin{equation}\eqlabel{WZWcdim}
\Delta_\lambda = \frac{\frac{1}{2}\left(\lambda,\lambda+2\rho\right)}{k+h^\vee},
\end{equation}
 where $\rho$ is the Weyl vector, i.e. $\rho=(1,1,\ldots,1)$ in the basis of fundamental weights, and $h^\vee$ is the dual Coxeter number, i.e. $h^\vee=N$ for $\mathfrak{su}(N)$.
Recall that $(\mu,\nu)=\sum_{i,j}\mu_i\nu_jG_{ij}$, where $G$ is the (symmetrized) inverse of the Cartan matrix. For $\mathfrak{su}(N)$ we have (see for example eq. $(7.15)$ of~\cite{FSbook})
\begin{equation}\eqlabel{cartinv}
G_{ij}=\frac{1}{N}\mathrm{min}(i,j)\left(N-\mathrm{max}(i,j)\right).
\end{equation}
A simple calculation shows that
\begin{equation}\eqlabel{Jpdim}
\Delta_{J^p} = \frac{p(N-p)k}{2N}.
\end{equation}
All invertible objects of odd order lie in $\mathrm{Pic}^\circ(\mathcal{C})$ (follows since every odd order element is a square, Remark 2.21 of~\cite{FRS3}), so it is enough to check whether even order elements satisfy~\equref{scineffcenter}.
First, note that if $k\equiv 0\mod N$, then $2\Delta_{J^p}\in\mathbb{Z}$, so all powers of $J$ lie in $\mathrm{Pic}^\circ(\mathcal{C})$, i.e. $\mathrm{Pic}^\circ(\mathcal{C})=Pic(\mathcal{C})\cong\mathbb{Z}/N\mathbb{Z}$.
More generally we have
\begin{proposition}\prlabel{SUNeffcenter}
(i) If $N$ is odd, or if $N$ and $k$ are both even, then $\mathrm{Pic}^\circ(\mathcal{C})=Pic(\mathcal{C})\cong\mathbb{Z}/N\mathbb{Z}$.
(ii) If $N=2p$ is even and $k$ is odd, then $\mathrm{Pic}^\circ\cong\mathbb{Z}/p\mathbb{Z}$, and $\mathrm{Pic}^\circ$ is generated by $J^2$. In particular, $\mathrm{Pic}^\circ(\mathcal{C})$ is closed under multiplication.
\end{proposition}
\begin{proof}
Assume $N$ is odd, and let $p\in\{1,\ldots,N-1\}$. If $\gcd(p,N)=1$ then $|J^p|=N$, which is odd and hence $J^p\in \mathrm{Pic}^\circ(\mathcal{C})$. If $\gcd(p,N)=q\neq 1$, then $q$ and $N/q$ are both odd, and we have $|J^p|=N/q$, hence $J^p\in \mathrm{Pic}^\circ(\mathcal{C})$. Assume instead $N$ and $k$ are both even, then $N\Delta_{J^p}=\frac{p(N-p)k}{2}\in\mathbb{Z}$ so $J^p\in \mathrm{Pic}^\circ(\mathcal{C})$ if $\gcd(p,N)=1$. If $\gcd(p,N)=q\geq 1$ then $|J^p|\Delta_{J^p}=\frac{N}{q}\frac{p(N-p)k}{2N}=\frac{p}{q}\left(N-p\right)\frac{k}{2}\in\mathbb{Z}$ and thus $J^p\in \mathrm{Pic}^\circ(\mathcal{C})$. This finishes the proof of (i). For part (ii), assume $N$ is even and $k$ is odd and let $\gcd(p,N)=q$. If $q=1$ then $|J^p|=N$, and $N\Delta_{J^p}=\frac{p(N-p)k}{2}$ For $p$ odd this is not an integer, while for $p$ even it is. If $q\neq 1$ and $p$ is odd, then also $p/q$ is odd. We then have $|J^p|\Delta_{J^p}=\frac{N}{q}\frac{p(N-p)k}{2N}=\frac{p}{q}\left(N-p\right)\frac{k}{2}\notin\mathbb{Z}$. Thus $J^p\notin \mathrm{Pic}^\circ(\mathcal{C})$ for $p$ odd. It remains to investigate the case $q\neq 1$ and $p$ even. If $N/q$ is odd we know $J^p\in \mathrm{Pic}^\circ(\mathcal{C})$. If instead $N/q$ is even, then $|J^p|\Delta_{J^p}=\frac{N}{q}\frac{p(N-p)k}{2N}=\frac{p}{2}\left(\frac{N-p}{q}\right)k\in\mathbb{Z}$. Thus $J^p\in \mathrm{Pic}^\circ(\mathcal{C})$ if $p$ is even.
\end{proof}
A Kreuzer-Schellekens bihomomorphism on the group $H\subset\mathrm{Pic}^\circ(\mathcal{C})$ satisfies $\Xi(g,g)=\theta_g$, so if $H$ is cyclic with generator $J$ the bihomomorphism is uniquely determined by $\Xi(J,J)=\theta_J$. Obviously we have $\Xi(J^a,J^b)=\theta_J^{ab}$ for any $a,b\in\mathbb{N}_0$.
In other words, for any cyclic subgroup $H\subset\mathrm{Pic}^\circ(\mathcal{C})$ there is a unique Schellekens algebra with support $H$.

\paragraph{$\mathbf{N=2}$}
This case was analysed in detail in~\cite{AF}, and we will very briefly review the results.
We are concerned with the categories $\mathcal{C}_{2,k}$ for $k\in\mathbb{N}$, where the isomorphism classes of simple objects are in bijection with the set $\{(0), (1), \ldots ,(k)\}$. Use the notation $(n)$, $n\in\{0,\ldots,k\}$ both for an isomorphism class, and for a representative simple object of this class. According to \prref{SUNeffcenter} we have 
\begin{equation*}
\mathrm{Pic}^\circ(\mathcal{C}_{2,k})=\left\{
	\begin{array}{ll}
	\{(0)\} & \text{if } k\not\in 2\mathbb{N}\\
	\{(0),(k)\} & \text{if } k \in 2\mathbb{N}
	\end{array}
	\right.
\end{equation*}
There are thus no non-trivial Schellekens algebras for odd $k$. We know more than this, however. The simple symmetric special Frobenius algebras in $\mathcal{C}_{2,k}$ have been classified up to Morita equivalence~\cite{KO, O}, with the result that for odd $k$ there is only the Morita class of the trivial algebra $\mathbf{1}=(0)$. The same holds for $k=2$; although there is a Schellekens algebra with support $\mathbb{Z}/2\mathbb{Z}$, this turns out to be Morita equivalent to $\mathbf{1}$. For every even $k$ larger than $2$ on the other hand, there are simple algebras not Morita equivalent to $\mathbf{1}$. The Schellekens algebra on the object $A=(0)\oplus (k)$ exists for every even $k$, and for $k\geq 4$ is not Morita equivalent to $\mathbf{1}=(0)$. In addition for $k=10, 16, 28$ there is a third Morita class of algebras. For $k=10$ an algebra can be constructed on the object $(0)\oplus (6)$, for $k=16$ there is an algebra on the object $(0)\oplus (8)\oplus (16)$, and for $k=28$ there is an algebra on the object $(0)\oplus (10)\oplus(18)\oplus(28)$, all of which belong to different Morita classes than the trivial algebra or the Schellekens algebra.This classification coincides precisely with the ADE classification of modular invariants for the same categories~\cite{CIZ}, and algebras not Morita equivalent to $\mathbf{1}$ give non-trivial modular invariants $Z(A)$~\cite{BE}. Hence the representations $V_g$ are reducible for even level $k\geq 4$.

We also showed in~\cite{AF} how one can use these techniques to determine dimensions of subrepresentations explicitly, and did this for the genus $1$ representations. When $k=4n$ there are subrepresentations of dimensions $n+1$ and $3n$ respectively, and when $k=4n+2$ there are subrepresentations of dimensions $n+1$ resp. $3n+2$. In addition, when $k=10$ the dimension $8$ part decompose in two subrepresentations of dimensions $3$ and $5$, and similar decompositions occur also when $k=16$ and $k=28$.

\paragraph{$\mathbf{N\geq 3}$}

\begin{proposition}\prlabel{levelN}
For level $k\equiv 0\mod N$, $N>2$, the Schellekens algebra $A$ in $\mathcal{C}_{N,k}$ with support $\mathbb{Z}/N\mathbb{Z}$ has non-trivial torus partition function $Z(A)$.
\end{proposition}
\begin{proof}
First note that if $N$ is even then also $k$ is even, so $\mathbb{Z}/N\mathbb{Z}$ lies in the effective center in all cases considered in the proposition. Since $\chi_{\mathbf{1}}(J)=1$ we have the following expression for $Z_{0,j}(A)$
\begin{equation}\eqlabel{pnfn1}
Z_{0,j}(A)=\frac{1}{N}\sum_{a,b=0}^{N-1}\Xi(J^a,J^b)\delta_{j,J^b}.
\end{equation}
Assume $N$ is odd, then $\Delta_J=\frac{k(N-1)}{2N}\in\mathbb{Z}$ since $k\propto N$ and $N-1$ is even. Thus for $N$ odd, $\Xi(J,J)=\theta_J=1$, and hence $\Xi=1$. It follows immediately that for any $b\in\{1,\ldots,N-1\}$, $Z_{0,J^b}(A)=1$.
Assume $N$ is even. Then $2\Delta_{J}\in\mathbb{Z}$, and it follows that the character $\Xi(\cdot,J^2)=1$ implying $Z_{0,J^2}(A)=1$ By assumption $N>2$, so $J^2\not\cong\mathbf{1}$ so $Z(A)$ is non-trivial.
\end{proof}
\begin{proposition}\prlabel{oddN}
If $N$ is odd, the Schellekens algebra $A$ with support $\mathbb{Z}/N\mathbb{Z}$ has non-trivial torus partition function $Z(A)$ for any level $k\in\mathbb{N}$.
\end{proposition}
\begin{proof}
If $k\equiv 0\mod N$ the statement follows from Proposition~\prref{levelN}, thus assume $N\not | k$. We have $\Xi(J,J^b)=\theta_j^b$, and $b\Delta_J=\frac{b(N-1)k}{2N}=\frac{bpk}{2p+1}$ if $N=2p+1$. If $\gcd(k,N)=q\neq 1$ then take $b=N/q$, and we have $N\Delta_J/q=p\frac{k}{q}\in\mathbb{Z}$, thus $\Xi(J,J^{N/q})=1$ so $\Xi(\cdot,J^{N/q})$ is the trivial character. The formula \equref{pnfn1} then implies that $Z_{0,J^{N/q}}=1$, and since $J^{N/q}\not\cong \mathbf{1}$ this shows that $Z(A)$ is non-trivial. If on the other hand $\gcd(k,N)=1$, then there is no $b\in\{1,\ldots,N-1\}$ such that $\Xi(\cdot,J^b)=1$, and therefore $Z_{0,j}(A)=\delta_{j,0}$. To see if $Z(A)$ is non-trivial we must then check coefficients $Z_{i,j}(A)$ where $i,j\neq 0$. Let $X=(1,0,\ldots,0)$ (which coincides with $J$ iff $k=1$), implying that $JX=(k-1,1,0,\ldots,0)$. It is straightforward to determine the conformal weight of $X$ and $JX$, with the results $\Delta_X=\frac{N^2-1}{2N(k+N)}$, $\Delta_{JX}=\frac{(k+1)N(N-2)+k(N-2)+k^2(N-1)-1}{2N(k+N)}$. Another simple calculation shows $\Delta_J+\Delta_X-\Delta_{JX}=\frac{1}{N}$, implying $$\chi_X(J)=\exp(2\pi i(\Delta_J+\Delta_X-\Delta_{JX}))=\exp(2\pi i/N).$$ Note that $\gcd(p,2p+1)=1$, so under the assumption $\gcd(k,N)=1$ we have $\gcd(pk,N)=1$ implying that $pk\mod N$ generate $\mathbb{Z}/N\mathbb{Z}$. There thus exists a $b\in\{1,\ldots,N-1\}$ such that $bpk\equiv 1\mod N$, and we get
$$\chi_X(J^a)\Xi(J^a,J^b)=\exp(2\pi i\frac{a}{N})\exp(-2\pi i\frac{abpk}{N})=\exp(2\pi i a\frac{1-bpk}{N})=1,\ \forall a\in\{1,\ldots,N-1\}.$$
With this choice of $b$ we then have
$$Z_{X,J^bX}=\frac{1}{N}\sum_{a,c=0}^{N-1}\chi_X(J^a)\Xi(J^a,J^c)\Delta_{J^bX,J^cX}=\frac{1}{N}\sum_{a=0}^{N-1}\chi_X(J^a)\Xi(J^a,J^b)=1.$$ Since for any $b\in\{1,\ldots,N-1\}$ we have $J^bX\not\cong X$, it follows that $Z(A)$ is non-trivial.
\end{proof}

\begin{proposition}\prlabel{evenNk}
Let $N$ and $k$ be even, and $N>2$. If $4\not | k$ or if $\gcd(N/2,k/2)\neq 1$, then the Schellekens algebra $A$ in $\mathcal{C}_{N,k}$ with support $\mathbb{Z}/N\mathbb{Z}$ has non-trivial torus partition function. If $4|k$ and $\gcd(N/2,k/2)=1$, then the Schellekens algebra with support $\langle J^2\rangle\cong\mathbb{Z}/\frac{N}{2}\mathbb{Z}$ has non-trivial torus partition function.
\end{proposition}
\begin{proof}
Set $N=2p$ and $k=2l$, then $\Delta_J=\frac{l(2p-1)}{2p}$. If $\gcd(l,p)=q\neq 1$, then $\frac{N}{q}\Delta_J=\frac{l}{q}(2p-1)\in\mathbb{Z}$ implying $\Xi(J,J^{N/q})=1$, hence $\Xi(\cdot,J^{N/q})$ is the trivial character. We thus have
$$Z_{0,J^{N/q}}(A)=\frac{1}{N}\sum_{a=0}^{N-1}\Xi(J^a,J^{N/q})=1,$$
and since $J^{N/q}\not\cong\mathbf{1}$, $Z(A)$ is non-trivial.

Assume instead $\gcd(l,p)=1$, and first consider the case $k\not\equiv 0\mod 4$. We then have $\gcd(l,N)=1$ (since $l$ odd). Since $\gcd(p,2p-1)=1$ we also have $\gcd(l(2p-1),N)=1$, and $l(2p-1)\mod N$ therefore generate $\mathbb{Z}/N\mathbb{Z}$. Let again $X=(1,0,\ldots,0)$. The same calculation as in the proof of~\prref{oddN} shows $\Delta_J+\Delta_X-\Delta_{JX}=\frac{1}{N}$, and that there exists $b\in\{1,\ldots,N-1\}$ such that $Z_{X,J^bX}(A)=1$. It therefore follows that $Z(A)$ is non-trivial.
Finally assume $k\equiv 0\mod 4$, and consider the algebra with support $\mathbb{Z}/p\mathbb{Z}$. Note that the corresponding Kreuzer-Schellekens bihomomorphism satisfies $\Xi(J^2,J^2)=\theta_{J^2}$. We have $$\chi_{J^{p-1}}(J^2)\Xi(J^2,J^2)=\exp(2\pi i[\Delta_{J^2}+\Delta_{J^{p-1}}-\Delta_{J^{p+1}}-\Delta_{J^2}])=\exp(2\pi i[\Delta_{J^{p-1}}-\Delta_{J^{p+1}}]),$$
and since $N=2p$, $$\Delta_{J^{p-1}}-\Delta_{J^{p+1}}=\frac{(p-1)(N-p+1)k-(p+1)(N-p-1)}{2N}=\frac{(p-1)(p+1)k-(p+1)(p-1)k}{2N}=0.$$ The homomorphism property then implies $\chi_{J^{p-1}}(J^{2a})\Xi(J^{2a},J^2)=1$ for every $a=1,\ldots,p-1$.
The formula $$Z_{J^{p-1},j}(A)=\frac{1}{p}\sum_{a,b=0}^{p-1}\chi_{J^{p-1}}(J^{2a})\Xi(J^{2a},J^{2b})\delta_{j,J^{2b}J^{p-1}}$$ then immediately implies $Z_{J^{p-1},J^{p+1}}(A)=1$. When $N=2p>2$, $J^{p-1}\not\cong J^{p+1}$, and it follows that $Z(A)$ is non-trivial. This finishes the proof.
\end{proof}

\begin{proposition}\prlabel{evenNoddk}
Let $N$ be even, $N>2$, and let $k$ be odd. Then the Schellekens algebra with support $\langle J^2\rangle\cong\mathbb{Z}/\frac{N}{2}\mathbb{Z}$ has non-trivial torus partition function.
\end{proposition}
\begin{proof}
Set $N=2p$. 
If $\gcd(k,p)=q\neq 1$ then $\Xi(J^2,J^{2p/q})=\theta_{J^2}^{p/q}=\exp(-2\pi i\frac{p}{q}\frac{2(N-2)k}{2N})=\exp(-2\pi i\frac{p}{q}\frac{(p-1)k}{p})=\exp(-2\pi i (p-1)\frac{k}{q})=1$. The formula $$Z_{0,j}(A)=\frac{1}{p}\sum_{a,b=0}^{p-1}\Xi(J^{2a},J^{2b})\delta_{j,J^{2b}}$$ then implies $Z_{0,J^{N/q}}(A)=1$. For $N>2$, $J^{N/q}\not\cong\mathbf{1}$ so $Z(A)$ is non-trivial. If $\gcd(k,p)=1$, consider instead $Z_{J^{p-1},j}$. One quickly verifies that $\chi_{J^{p-1}}(J^{2a})\Xi(J^{2a},J^{2b})=1$ for $a,b\in\{1,\ldots,p-1\}$, and it follows that $Z_{J^{p-1},J^{p+1}}(A)=1$. Since $J^{p-1}\not\cong J^{p+1}$ if $N>2$, $Z(A)$ is non-trivial.
\end{proof}

The results~\prref{oddN}, 3.4, and 3.5 imply the following corollary to~\thref{red}
\begin{theorem}\thlabel{redSUN}
The quantum representations $V_g$ constructed from $\mathbf{tqft}_\mathcal{C}$ where $\mathcal{C}$ is the category of integrable highest weight representations of $\widehat{\mathfrak{su}}(N)$ for any level $k\in\mathbb{N}$ and $N>2$ are reducible for all genus $g\geq 1$.
\end{theorem}

\begin{remark}
When $N=3$ one easily compares the torus partition functions of the Schellekens algebras with the complete list of modular invariants from~\cite{G}, with the expected result. When $k\geq 3$, the Schellekens algebra gives the D-type invariant multiplied with the charge conjugation matrix $C$. For $k=1, 2$, the Schellekens algebra gives $Z(A)=C$, which is still non-trivial.
\end{remark}

\end{document}